\documentclass[11pt]{amsart}
\usepackage{amsmath,amsfonts,amsthm,amssymb,amscd}
\usepackage[latin1]{inputenc}
\usepackage{appendix}
\usepackage{latexsym}
\usepackage{graphicx}
\usepackage{enumerate}

\def\M{\mathbb{N}_0}
\def\N{\mathbb{N}}
\def\R{\mathbb{R}}
\def\G{\text{GDWN}}
 
\newcommand\T{\rule{0pt}{2.6ex}}

\theoremstyle{definition}
\newtheorem{Thm}{Theorem}[section]
\newtheorem{Prop}[Thm]{Proposition}
\newtheorem{Lemma}[Thm]{Lemma}
\newtheorem{Cor}[Thm]{Corollary}
\newtheorem{Def}{Definition}

\newtheorem{Conj}{Conjecture}
\newtheorem{Ex}{Example}
\newtheorem{Rem}{Remark}

\title{A generalised diagonal wythoff nim}
\author{Urban Larsson}
\email{urban.larsson@chalmers.se} 
\address{Mathematical Sciences,
Chalmers University Of Technology and G\"oteborg University,
G\"oteborg, Sweden}
\date{\today }
\begin{document}
\begin{abstract}
In this paper we study a family of 2-pile Take Away games, that we denote by 
Generalized Diagonal Wythoff Nim (GDWN).
The story begins with 2-pile Nim whose sets 
of options and $P$-positions are $\{\{0,t\}\mid t\in \N\}$ and 
$\{(t,t)\mid t\in \M \}$ respectively. 
If we to 2-pile Nim adjoin the main-\emph{diagonal} 
$\{(t,t)\mid t\in \N\}$ as options, the 
new game is Wythoff Nim. 
It is well-known that the $P$-positions of this game 
lie on two 'beams' originating at the 
origin with slopes 
$\Phi= \frac{1+\sqrt{5}}{2}>1$ and $\frac{1}{\Phi } < 1$. 
Hence one may think of this as if, in the process of 
going from Nim to Wythoff Nim, the 
set of $P$-positions has \emph{split} and landed some distance off 
the main diagonal. 
This geometrical observation has motivated us to
 ask the following intuitive question. Does this splitting 
of the set of $P$-positions continue in some meaningful way 
if we, to the game of Wythoff Nim, adjoin some new \emph{generalized diagonal} 
move, that is a move of the form $\{pt, qt\}$, where $0 < p < q$ are fixed 
positive integers and $t > 0$? 
Does the answer perhaps depend on the specific values of $p$ and $q$? 
We state three conjectures of which the weakest form is: 
$\lim_{t\in \N}\frac{b_t}{a_t}$
exists, and equals $\Phi$, if and only if $(p, q)$ 
is a certain \emph{non-splitting pair}, and where $\{\{a_t, b_t\}\}$ represents 
the set of $P$-positions of the new game. 
Then we prove this conjecture for the special 
case $(p,q) = (1,2)$ (a \emph{splitting pair}).
We prove the other direction whenever $q / p < \Phi$.  
In the Appendix, a variety of experimental data is included, aiming to point 
out some directions for future work on GDWN games.

\end{abstract}
\maketitle
\vskip 30pt
\section{Introduction}
In this paper we analyze generalizations of the \emph{impartial} 
(see \cite{Con76, Lar09}) combinatorial games of 2-pile Nim \cite{Bou02} 
and Wythoff Nim, \cite{Wyt07, Fra82, FrOz98, HeLa06, Lar09, Lar, Lar2}. 
As usual, we let $\N$ denote the positive integers, $\M$ the non-negative 
integers and $\R$ the real numbers.

The options of 2-pile Nim are of the form $(x, y + t)\rightarrow (x, y)$ 
or $(x + t, y)\rightarrow (x, y)$, $t\in \N$, $x,y\in \M$, and the 
$P$-positions are of the form $(t,t)$, $t\in \M$. 
These positions may be thought of, geometrically, as 
one singular infinite North-East $P$-\emph{beam},
originating at the origin. 
In the game of Wythoff Nim a player may move as in Nim but also 
$(x+t,y+t)\rightarrow (x,y)$, $t\in \N$. 
For this game, geometrically, the 
singular Nim-beam of $P$-positions has \emph{split} into precisely 
two distinct beams, still originating at the origin, 
one leaning towards 'North' and the other towards 'East'. 
Let $\Phi = \frac{\sqrt{5}+1}{2}$ denote the Golden ratio. 
It is well-known that a position of this game is a $P$-position if and only if
it is of the form 
$(\lfloor \Phi t\rfloor ,\lfloor(\Phi + 1)t\rfloor )$ or
$(\lfloor (\Phi +1)t\rfloor , \lfloor \Phi t \rfloor )$, so that the 
new $P$-beams have slopes $\Phi$ and $\Phi - 1$ respectively. 

This geometrical observation has motivated us to ask the following 
intuitive question. 
Does this 'splitting behavior', going from Nim to Wythoff Nim, continue 
in some meaningful way 
if we, to the game of Wythoff Nim, adjoin some \emph{generalized diagonal} 
move, that is a move of the form $\{pt,qt\}$, where $0 < p < q$ are fixed 
positive integers and $t>0$, 
and then play the new game with both the old and the new moves?
Does the answer perhaps depend on the specific values of $p$ and $q$? 
Here we only study symmetric game rules so that in the coming, for a 
specific game, $(x,y)$ 
is a $P$-position if and only if $(y,x)$ is also. We indicate this by 
rather denoting such positions by $\{x, y\}$.

In Section 2, we state three conjectures for our family of new games 
that we denote Generalized Diagonal Wythoff Nim, \G.
The weakest form of the conjectures is: 
$$\lim_{t\rightarrow \infty}\frac{b_t}{a_t}$$ 
exists, and equals $\Phi$, if and only if $(p, q)$ 
is a non-splitting pair. Here $\{\{a_t,b_t\}\}$ represents the set of 
$P$-positions of the new game (with $(a_t)$ increasing) and 
by a \emph{splitting pair} we mean a pair of integers of the form 
$(\lfloor \Phi t\rfloor,  \lfloor (\Phi +1) t\rfloor) \text{ or } 
(\lceil \Phi t\rceil, \lceil (\Phi+1) t\rceil)$, $t>0$. 
Then we prove this conjecture for the special 
case $(p,q) = (1,2)$ (a splitting pair) which is the main result of Section 3.

In Section 4 we prove the other direction of this conjecture for a 
whole subfamily of games, namely whenever $q/p < \Phi$.

To give a hint of the direction of this work, we begin by presenting 
a table of $P$-positions and three figures. More of this kind may 
be found in the Appendix.
\begin{table}[ht!]
\begin{small}
\begin{center}
\setlength{\tabcolsep}{3pt}{
  \begin{tabular}
{| l || c|c|c|c| c| c| c| c| c| c| c| c| c| c| c| c| c| c| c|}
    \hline
    $b_n$\T   &0&\bf{3}&\bf{6}&5&10&14&17&\bf{25}&\bf{28}&18&\bf{35}&23&31&29&\bf{48}&32&\bf{55}&37&40\\
    $a_n$\T   &0&1&2&4&7 &8 &9 &11&12&13&15&16&19&20&21&22&24&26&27\\\hline
 $\delta_n$\T &0&2&4&1&3 &6 &8 &14&16&5 &20&7 &12&9 &27&10&31&11&13\\
$\gamma_n$\T  &0&1&2&-3&-4&-2&-1&3&4&-8&5&-9&-7&-11&6&-12&7&-15&-14\\
$\eta_n$\T  &0&5&10&6&13&20&25&39&44&23&55&30&43&38&75&42&86&48&53\\\hline
    $n$ \T    &0&1&2&3&4 &5 &6 &7 & 8&9 &10&11&12&13&14&15&16&17&18\\
    \hline
  \end{tabular}}
\end{center}
\end{small}
\caption{Here $\{a_n,b_n\}$ represents a $P$-position of $(1,2)\G$ 
for $n\in [0,18]$. 
Further, $\delta_n=b_n-a_n$, $\gamma_n=b_n-2a_n, \eta_n=2b_n-a_n$.}
\end{table}

\begin{figure}[ht!]
\begin{center}
	  {\includegraphics[width=0.5\textwidth]{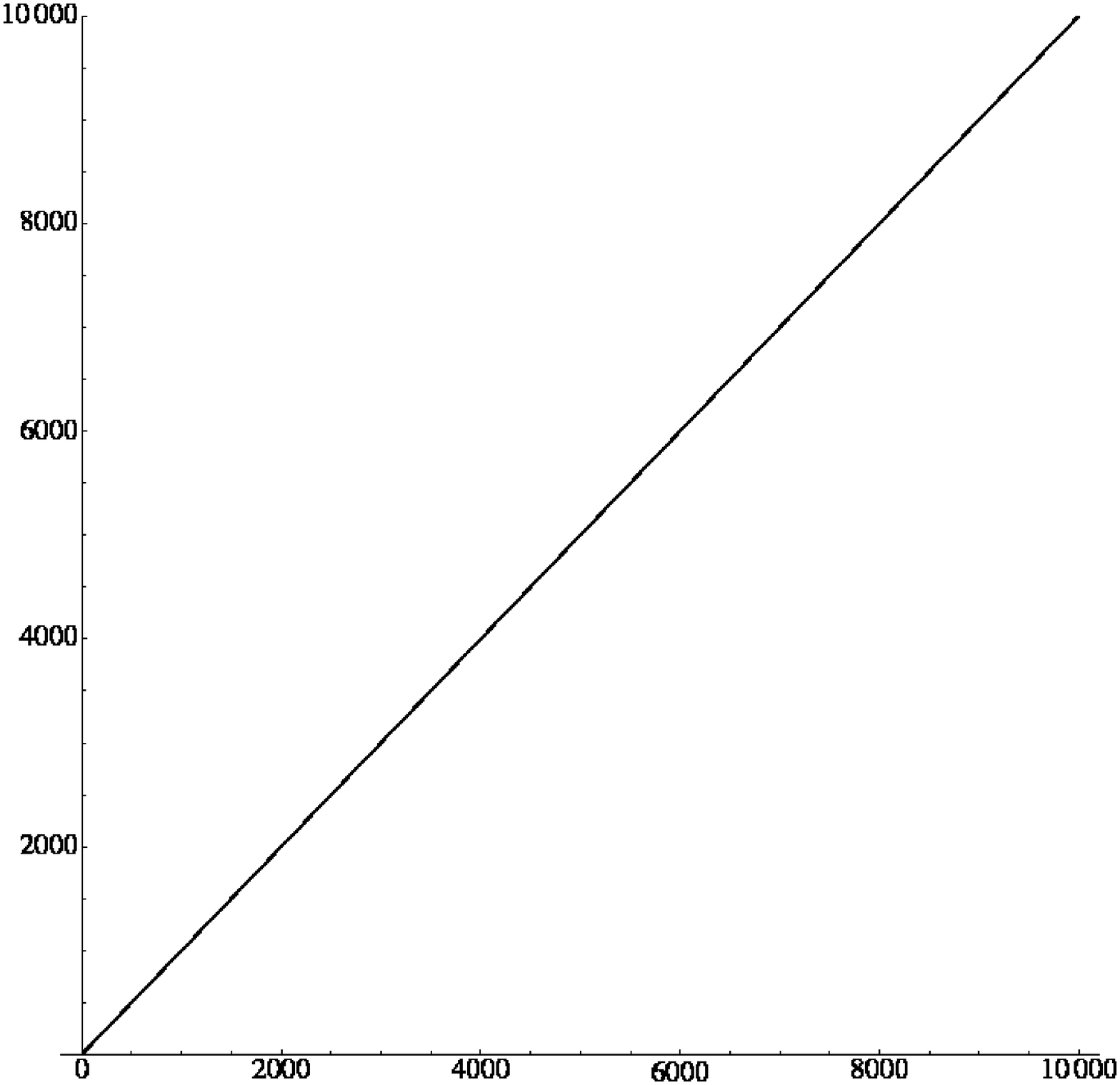}}  
\end{center}\caption{The $P$-positions of 2-pile Nim, the single ``Nim beam'' 
of slope 1.}
\begin{center}
	  {\includegraphics[width=0.5\textwidth]{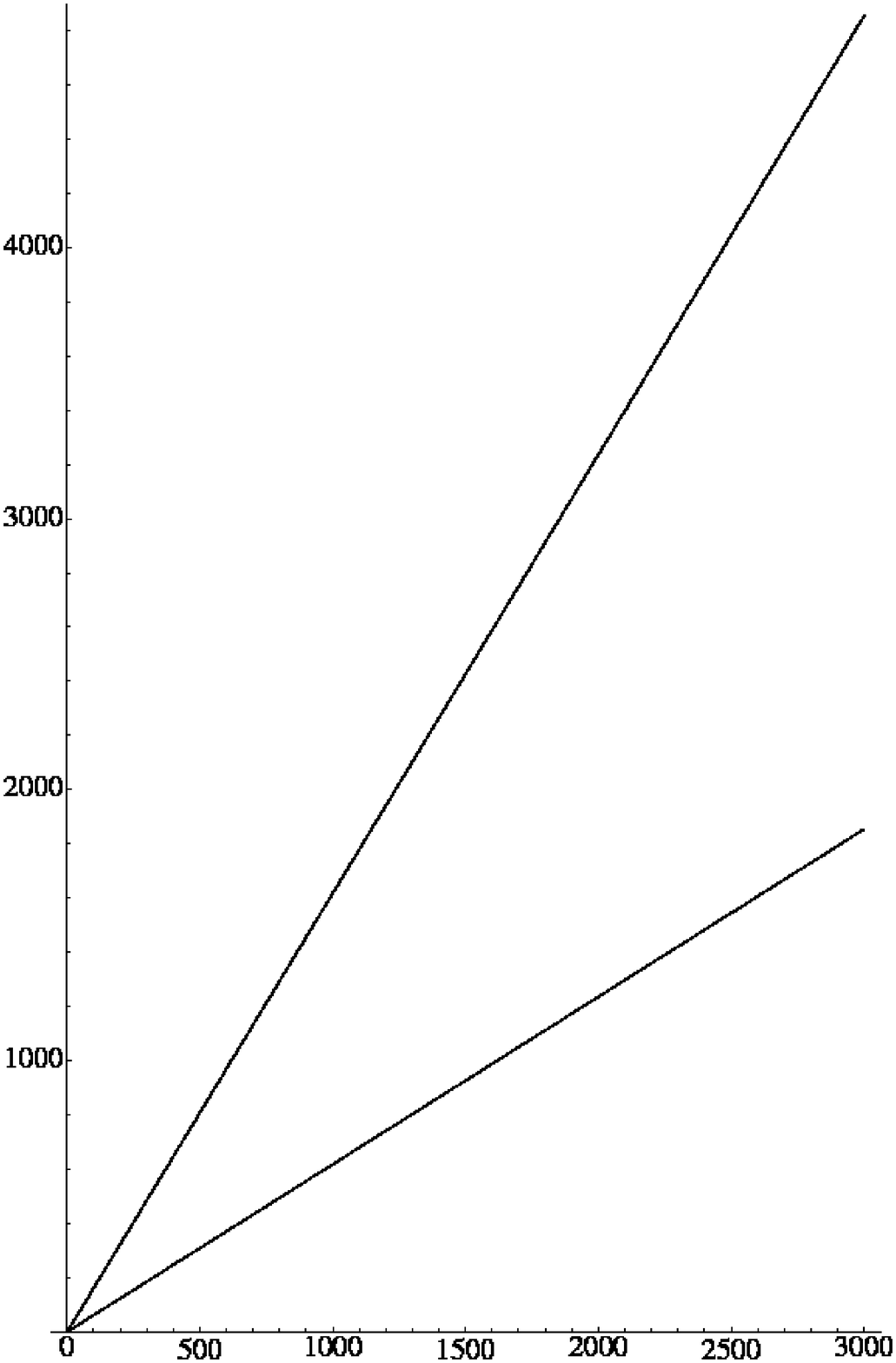}}  
\end{center}\caption{The $P$-positions of Wythoff Nim. This illustrates 
what we in Example 3 call the fundamental $1$-split. These $P$-beams have 
slopes $\frac{1+\sqrt{5}}{2}$ and $\frac{2}{1+\sqrt{5}}$ respectively}
\end{figure}
\clearpage
\begin{figure}[ht!]
\begin{center}
	  {\includegraphics[width=0.65\textwidth]{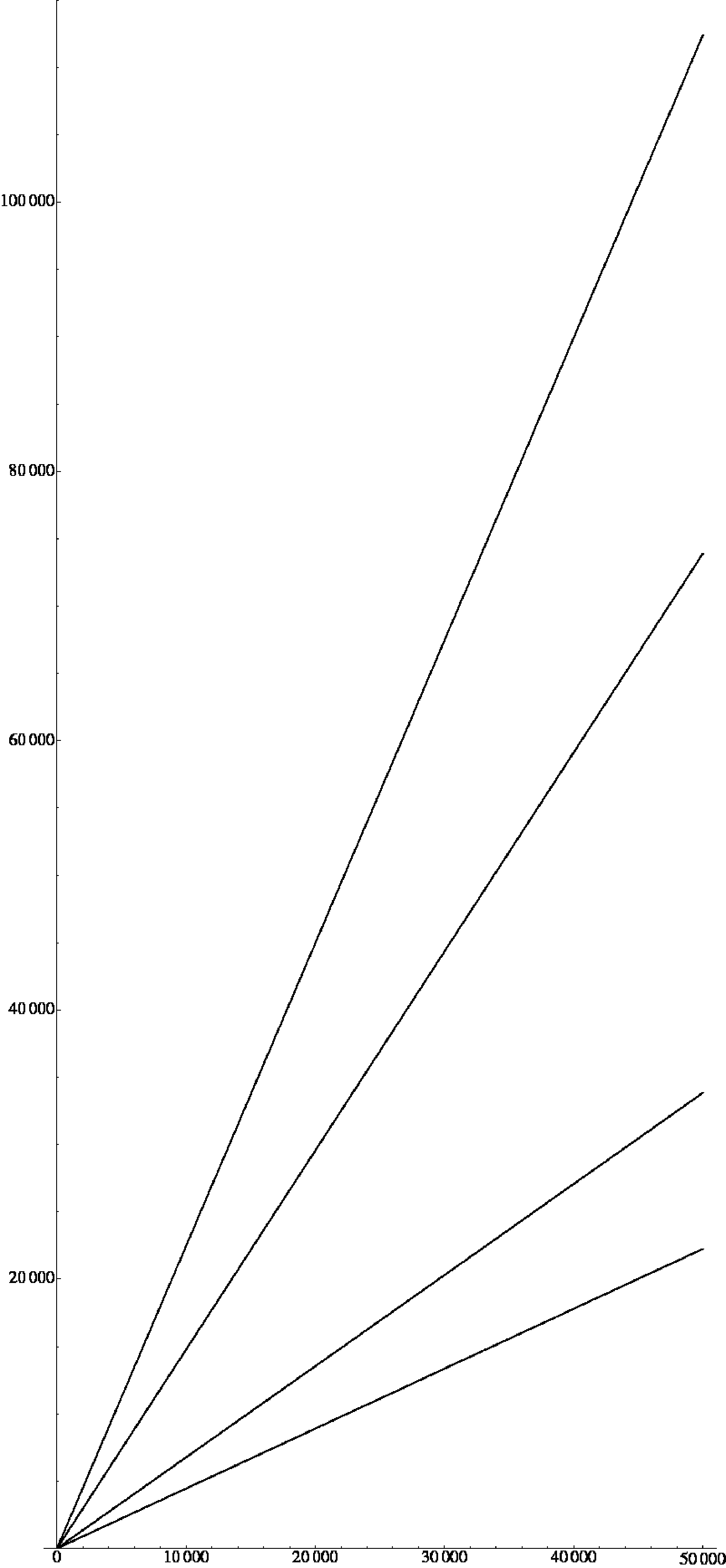}}  
\end{center}\caption{The $P$-positions $\{a_n, b_n\}$ of $(1,2)\G$ 
and $0\le n \le 50000$. Our computations seem to suggest that the 
slopes of the upper two $P$-beams 
are $2.247\ldots$ and $1.478\ldots$ respectively, see also Figure \ref{A12}.}
\end{figure}
\clearpage
\section{Sequences, games and conjectures}
 We begin with a general definition 
of the games and sequences explored in this paper. This definition is 
a straightforward generalization of 2-pile Nim and Wythoff Nim (see also 
Example 1 and 2 below).

\begin{Def}\label{Def:1}
Let $\mathcal{Q}_k$, $k\in \M$, denote the family of all sets of pairs 
of integers of the form 
$\{(p_i, q_i)\mid i\in \{0,1,\ldots ,k\}, p_i\in \M, q_i\in \N, 
p_i\le q_i, (p_i,q_i)\ne t(p_j,q_j) \text{ if } i\ne j, t\in \N \}$. 
Let $\mathcal{Q} := \cup_{k\in\M}\mathcal{Q}_k$.
\begin{itemize}
\item [(i)] Let $k\in \M$, $(x, y)\in \M\times\M$ and $Q\in \mathcal{Q}_k$. 
Then $(x, y)\rightarrow (x - m, y - n)$ is 
a legal move, of a game that we denote by $Q\G$, if 
$x - m\ge 0$ and $y - n\ge 0$ and if, 
for some $t\in \N$ and some $i\in \{0,1,\ldots ,k\}$, either 
\begin{align*}
m = tp_i \text{ and } n = tq_i \intertext{or} m = tq_i \text{ and } n = tp_i.
\end{align*}
The whole family of such games is simply denoted by \G. 
\item [(ii)] Define a function $\pi:=\pi_Q:\M \rightarrow \M$ 
recursively as: $\pi(n)$ is the least non-negative number distinct from 
 $$\frac{q_i\pi(j) + p_i(n - j)}{q_i},$$ and from
 $$\frac{p_i\pi(j) + q_i(n - j)}{p_i}$$ whenever $q_i\ne p_i > 0$,
for all $i\in \{0,1,\ldots ,k\}$ and for all $j\in \{0,1,\ldots ,n-1\}$. 
\end{itemize}
\end{Def}

\begin{Rem}
There is a reason for start indexing the $(p_i,q_i)$:s with zero 
rather than one. For the purpose of this paper and as we explain 
in the paragraph just after Example 2, 
we will insist that $(p_0, q_0) := (0, 1)$. (This may be put in contrast with 
\cite{FHL}, where we study games void of Nim-type moves). With respect to the
conjecture in Remark \ref{Rem:4} and possible future work, 
it will be convenient that $k$ counts the number of 'non-Nim type' 
diagonals adjoined. Here, almost exclusively, we 
will restrict our attention to the case $k=2$.
\end{Rem}

\begin{Def}\label{Def:2}
Let $\tau:\M\rightarrow \M$. Define $U = U_{\tau}=(u_i)_{i\in \M}$ as 
the increasing sequence of all $i$ such that $\tau(i)\ge i$ and similarly, 
with $L = L_{\tau} = (l_i)_{i\in \M}$ the increasing sequence of all $i$ such 
that $\tau (i) < i$. For a fixed $Q\in \mathcal{Q}$, put 
$U_Q = U_{\pi_Q}$ and $L_Q = L_{\pi_Q}$. 
\end{Def}

It follows immediately that $U\cup L = \M$ and $U\cap L = \emptyset$. 

\begin{Ex} \label{Ex:1}
As we have hinted, 2-pile Nim and Wythoff Nim are special cases of \G.
\begin{itemize}
\item The game $(0,1)\G$ is 2-pile Nim. For all $n\in \M$, $\pi_{(0,1)}(n) = n$.
\item The game $(0,1)(1,1)\G$ is Wythoff Nim. For all $n\in U_\pi $ 
we have $\pi(n) = \lfloor \Phi n \rfloor$. Otherwise 
$\pi(\lfloor \Phi n\rfloor) = n$. (See also \cite{HeLa06}.)
\end{itemize}
\end{Ex}

\begin{Def}\label{Def:3}
Suppose that $\tau = \pi_Q$, for some $Q$ as in Definition \ref{Def:1}, 
and define $a = a(Q) = (a_i)$ and $b = b(Q) = (b_i)$ by, for all $i$, 
$a_i = u_i$ and $b_i = \pi_Q(u_i)$,
\end{Def}

Then $a_0 = b_0 = 0$, $a_1=1$, $a$ and $b\cap \N$ are \emph{complementary}, 
that is $a\cup b = \M$ and $a\cap b = \{0\}$, $a$ is increasing, but, 
in general, $b$ is not. Is it true that $b(Q) = L_Q$, that is that $b$ 
is increasing, if and only if 
$\{(n, \pi_Q(n))\mid n\in \M\} = \mathcal{P}(\text{Wythoff Nim})$?

\begin{Def}
Suppose $G$ is an \emph{impartial} game. Then $G$ 
is $P$ if none of the options of $G$ is $P$. Otherwise, $G$ is $N$.
\end{Def}

An immediate consequence of this definition is that the next player to 
move wins if and only if $G$ is $N$.
We denote with $\mathcal{P}(G)$ the complete set of $P$-positions of $G$.

\begin{Thm}\label{Thm:2.1}
Fix a $Q \in \mathcal{Q}$. Then 
\begin{itemize}
\item [(i)] $\pi_Q$ is an involution of $\M$, that is, for all $i$, 
$\pi(i) = \pi^{-1}(i)$.
\item [(ii)] $\mathcal{P}(Q\G) = \{(i, \pi(i))\mid i\in \M \} = 
\{\{a_i,b_i\}\mid i\in \M \}$.
\end{itemize}
\end{Thm}

\noindent{\it Proof.} This is immediate by definition. \hfill $\Box$

\begin{Ex}\label{Ex:2} Some games are particularly easy to analyze. The 
first two items are the same as in Example \ref{Ex:1}. 
\begin{itemize}
\item The set of $P$-positions of 2-pile Nim is 
$\{(n,n)\mid n\in \M \} = \{(n,\pi(n))\mid n\in \M \}$.
\item The $P$-positions of Wythoff Nim are usually represented as 
all pairs of the form $\{A_n, B_n\}$, where $A_n:= \lfloor n\Phi\rfloor$ 
and $B_n := \lfloor n\Phi^2 \rfloor $, $n \in \M$. 
\item $\mathcal{P}((r,s)\G)=\{\{m, n\}\mid 0\le m,n < s\}\cup 
\{\{n,m\}\mid n < r \}$ if $r > 0$ and 
$\mathcal{P}((0,s)\G)=\{(sn+i,sn+j)\mid 0\le i,j < s, n\in \M\}.$
\item $\mathcal{P}((0,1)(r,s)\G) = \mathcal{P}($2-pile Nim$)$ 
if and only if $r\ne s$. For the case $r = s > 1$, on the one hand, 
the situation seems more complicated, see \cite{DuGr}, although probably 
still tractable, as discussed in \cite{FrPe}.
\item On the other hand, the solution of the games $(0,s)(s,s)\G$ may 
be represented via Beatty sequences, namely, 
$$\mathcal{P}((0,s)(s,s)\G) = \{\{sA_n+i, sB_n+j\} 
\mid 0\le i,j < s, n\in \M \}.
\footnote{We were not able to locate 
this game and its solution in the literature. We omit the proof, since 
it only uses an elementary inductive argument.}$$ 
\end{itemize}
\end{Ex}

From now onwards, we will only consider extensions of Wythoff Nim, 
that is games where the moves of Wythoff Nim is a subset of all legal 
moves for the new game. 
Hence, for the rest of this paper, let $(p_0,q_0)=(0,1)$ and $(p_1,q_1)=(1,1)$. 
In fact, except in the Appendix and in Remark \ref{Rem:3} and \ref{Rem:4}, 
we restrict our attention to games where $Q = \{(0,1),(1,1),(p,q)\}$ 
for some $p, q\in \N$ with $p < q$. Hence, to simplify notation, 
let $(p,q)\G$ denote $Q\G$ and let $\pi_{p,q}$ denote $\pi_{(0,1)(1,1)(p,q)}$. 
Then, with notation as in Example \ref{Ex:2}, 
we call $(A_n, B_n), n\in \N$, the \emph{Wythoff pairs}. 
The \emph{dual Wythoff pairs} are all pairs of the form  
$(\lceil n\Phi\rceil, \lceil n\Phi²\rceil ) = (A_n + 1, B_n + 1)$, $n\in \N.$ 

\begin{Rem} See also for example \cite{HeLa06} for an introduction to 
the closely related pair of arrays, the \emph{Wythoff Array} 
(the Wythoff pairs)
and the \emph{Dual of the Wythoff Array} (the Dual Wythoff pairs). In 
the Appendix we have given the first few entries of these arrays.
\end{Rem}

\begin{Def}\label{Def:4}
Let $p,q\in \N$. Then $(p,q)$ is a \emph{splitting pair} if it is 
a Wythoff pair or a Dual Wythoff pair.
\end{Def}
  
We will now present three conjectures and then prove some special 
formulations of the first.

\begin{Conj}\label{Conj:1}
Fix $p,q\in \N$, $a = a(p,q)$ and $b = b(p,q)$. Then the limit 
$$\lim_{n\in \N}\frac{b_n}{a_n}$$ 
exists and equals the Golden ratio, $\Phi$, if and only if $(p,q)$ 
is a non-splitting pair.
\end{Conj}

In light of 
a great deal of experimental data we may strengthen Conjecture \ref{Conj:1}. 
To this purpose we need to extend our terminology.
\begin{Def}
Let $\mu\in \R $, $\mu>0$. A sequence of pairs of positive 
integers $X=((x_i,y_i))_{i\in \N}$ with $x_i\le y_i$ $\mu$-\emph{splits} if 
there is an $\alpha \in \R$ such that,
\begin{itemize}
\item there are at most finitely many $i$:s 
such that $y_i/x_i\in [\alpha, \alpha +\mu )$,
\item there are infinitely many $i$:s such that 
$y_i/x_i\in [\alpha +\mu,\infty )$,
\item there are infinitely many $i$:s such 
that  $y_i/x_i\in (0, \alpha )$.
\end{itemize}
We say that $((x_i,y_i))_{i\in \N}$ \emph{splits} if there is a $\mu$ such 
that $((x_i,y_i))_{i\in \N}$ $\mu$-splits. 
If $X$ splits we may take $\xi \in [\alpha,\alpha +\mu)$ and define 
 complementary sequences $(l_i)$ and $(u_i)$ such that for all $i$
$$y_{l_i}/x_{l_i}\in (0, \alpha +\xi )$$ and 
$$y_{u_i}/x_{u_i}\in [\alpha +\xi ,\infty ).$$
\end{Def}
The most 'prominent splitting sequence' of the form $\mathcal{P}(Q)$, 
$Q\in \mathcal{Q}$, is the following example.
\begin{Ex}[The fundamental splitting sequence]\label{Ex:3}
Clearly $((i, \pi_{(0,1)}(i)))$ does not 
split, but $((i,\pi_{(0,1)(1,1)}(i)))$ does. Indeed, 
the latter $1$-splits (with $\alpha =\Phi - 1$). Take $\xi = 1$. Then, 
for all $i$, $b_i = l_i \in L$ and $a_i = u_i \in U$, 
where $Q = \{(0,1),(1,1)\}$ and where $a$ and $b$ are as 
in Definition \ref{Def:3}. 
\end{Ex}
Since a sequence of pairs of integers can split 'once', 
it is not unreasonable to think it could potentially 'split twice'.
\begin{Def}
Suppose that $((x_i,y_i))_{i\in \N}$ \emph{splits} and that $\alpha$ and $\mu$ 
is chosen so that $\mu$ is largest possible. Then, if there 
is an interval $[\beta,\beta +\mu']$ such that 
$[\beta,\beta +\mu']\cap [\alpha, \alpha +\mu )=\emptyset $ and such that 
\begin{itemize}
\item there are at most finitely many $i$:s 
such that $y_i/x_i\in [\beta, \beta +\mu' )$, and either
\item $\beta > \alpha +\mu$ and
\begin{itemize}
\item there are infinitely many $i$:s such that 
$y_i/x_i\in [\beta +\mu',\infty )$,
\item there are infinitely many $i$:s such 
that  $y_i/x_i\in [\alpha +\mu, \beta )$, 
\end{itemize}
\item or $\beta +\mu'< \alpha $, and
\begin{itemize}
\item there are infinitely many $i$:s such that 
$y_i/x_i\in [\beta +\mu',\alpha )$,
\item there are infinitely many $i$:s such 
that  $y_i/x_i\in [1, \beta )$, 
\end{itemize}
\end{itemize}
then we say that $((x_i,y_i))_{i\in \N}$ \emph{splits twice}. If 
$((x_i,y_i))_{i\in \N}$ splits, but does not split twice, we say that 
$((x_i,y_i))_{i\in \N}$ splits (precisely) \emph{once}. 
\end{Def}
\begin{Conj}\label{Conj:2}
Fix $p,q\in \N$ and define $\pi=\pi_{p,q}$, $a = a(p,q)$ and $b = b(p,q)$ 
as before. Then 
\begin{itemize}
\item [(i)] $((n,\pi(n)))_{n\in U} = ((a_n,b_n))_{n\in \M}$ splits if and only 
if $(p,q)$ is a splitting pair. 
\item  [(ii)] If $(p,q)$ is a splitting pair, then 
$((n, \pi(n)))_{n\in U}$ splits precisely once.
\end{itemize}
\end{Conj}

\begin{Rem}\label{Rem:3}
Let $k\in \N$, $k\ge 3$. Suppose that $((x_i,y_i))_{i\in \N}$ splits twice. 
By elaborating on the above definitions we may define conditions for 
$((x_i,y_i))_{i\in \N}$, a \emph{$k$-fold split}. And indeed, the next remark 
is supported by numerous computer simulations. A rigorous treatment 
of this and the next remark is left for future research. 
\end{Rem}

\begin{Rem}\label{Rem:4}
Fix $k\in \M$ and $Q = ((p_i,q_i))_{i\in \{0,1,\ldots ,k\}}$.
We conjecture that $((n,\pi_Q(n)))_{n\in \N}$ is a $C$-fold split with $C\le k$. 
Question: Is there, for each 
$k\in \M$, a set $Q\in \mathcal{Q}$ such that $\mathcal{P}(Q\G)$ is a $k$-fold
split?
\end{Rem}
By our simulations, the next more precise form of these conjectures is 
only applicable for certain values of $p$ and $q$, see Figure A12 and A14 
in the Appendix.

\begin{Conj}\label{Conj:3}
Fix $(p,q)=(1,2)$ or $(p,q)=(2,3)$. Then there is a pair of increasing 
complementary sequences $(l_i)$ and $(u_i)$ such that both 
$\eta = \lim_{i\rightarrow \infty}\frac{b_{l_i}}{a_{l_i}}$ and 
$\gamma = \lim_{i\rightarrow \infty}\frac{b_{u_i}}{a_{u_i}}$ 
exist with real $1< \eta < \Phi < \gamma \le 3$.
\end{Conj}
The upper bound on $\gamma$ is easy to verify (see Corollary 3.6).

Let us recall the first few $P$-positions of Wythoff Nim:
$$(0,0),(1,2),(3,5),(4,7),(6,10),(8,13),(9,15),(11,18),(12,20),(14,23),\ldots $$
In the Appendix we give the first few $P$-positions of $(p,q)\G$ 
for $$(p,q)=(1,2),(2,3),(2,4),(4,6) \text{ and }(4,7).$$ 
Notice that $(1,2)$ and $(4,7)$ are Wythoff pairs, $(2,3)$ and $(4,6)$ 
are dual Wythoff pairs whereas $(2,4)$ is neither. 
See also the Appendix for several plots of the ratios $b_i/a_i$ 
for different $p$ and $q$.

\section{A resolution of Conjecture \ref{Conj:1} for $(1,2)\G$}
Suppose $f,g: \M \rightarrow \R$. In this section we use the notation 
$f(N)\ll g(N)$ if $f(N) < g(N)$ for all sufficiently large $N$. 
(And analogously for $\gg$), where the term sufficiently large is explained 
by each surrounding context. We will have use for a simple but general lemma. 
\begin{Lemma}\label{Lemma:3.1}
Let $\tau:\M \rightarrow \M$ be an involution of the non-negative integers, 
that is, for all $i\in \M$, $\tau(\tau(i)) = i$. Then, for all $n\in \M$, 
 $$T:=\#\{i \mid i\in U_{\tau}\cap \{0,1,\ldots ,n\} \} \ge \frac{n+1}{2}.$$
\end{Lemma}
\noindent{\bf Proof.} 
Suppose on the contrary that $T < \frac{n+1}{2}$ for some $n\in \M$. 
Then $\frac{n+1}{2}\le \# \{i\mid i\in L_\tau \cap \{0,1,\ldots ,n\}\}$. Since 
$\tau(\tau(i))=i\le n$ and $i\in L_\tau$ gives $i>\tau(i)$, we get 
$\tau(i)\in U_\tau$ with $\tau(i)\le n$. Hence 
$$\frac{n+1}{2}
\le \#\{\tau (i) \mid \tau(i)\in U_\tau\cap \{0,1,\ldots ,n\} \}\le 
T < \frac{n+1}{2},$$ a contradiction. 
$\hfill \Box$\\
\begin{Lemma}\label{Lemma:3.2} Fix a $Q\in \mathcal{Q}$ and an $N\in \N$. 
Then $$\#\{i\mid a_i\le N\} > \frac{N}{2}.$$
\end{Lemma}
\noindent{\bf Proof.} 
This follows by Lemma \ref{Lemma:3.1} and Theorem \ref{Thm:2.1} (i) since 
$$\#\{i\mid a_i\le N\} = \#\{i\le N\mid i\in U_Q\}.$$ 
$\hfill \Box$\\
Some variation of the next Lemma has been studied before 
(see for example \cite{FrKr, HeLa06}). It is quite general, alas not as 
general as Lemma \ref{Lemma:3.1}. For our purpose we note that it holds for any 
$\tau=\tau_Q,$ where $\{(0,1),(1,1)\}\subset Q$.

\begin{Lemma}\label{Lemma:3.3}
Let $\tau$ be as in Lemma \ref{Lemma:3.1}. Suppose that, for all $i,j \in \M$, 
$\tau(i) - i = \tau(j) - j$ implies $i=j$.
Then the set $$\left\{i\in U_\tau \mid \frac{\tau(i)}{i}\ge \Phi \right\}$$ 
is infinite. 
\end{Lemma}
\noindent{\bf Proof.} 
Let $C\in \R$ with $C > 1$ and define the set
$$S = S(\tau, C) := \left\{i\in U_\tau \mid \frac{\tau(i)}{i}< C \right\}.$$
Suppose that $S$ contains all but finitely many elements of $U_\tau$. 
We have to show that $C > \Phi$. Put $c := C - 1 > 0$. 
Let $N\in \N$ be sufficiently large, 
by which we mean: For all $i\le N$ we have $\tau (i)\le CN$. 
(By the finiteness of $U \setminus S$ this is certainly possible.) 
Define $x\in \N$ by 
$u_x\le N < u_{x+1}$, where as before $(u_i)_{i\in \N}:=U_\tau $. 
Denote with $\delta_i = \tau(u_i) - u_i$. Since the $\delta_i$:s are distinct 
we must have $\max\{\delta_i\mid i\le x\}\ge x$. Then there exists 
an $i\le x$ such that $u_i + x\le \tau(u_i)$. But, by assumption, 
for a sufficiently large $N$ this implies 
$$1+\frac{x}{u_i}\le \frac{\tau(u_i)}{u_i} < C$$ and so 
\begin{align}\label{cN} 
x\ll (C-1)u_i\le cN.
\end{align}
On the other hand, by Lemma \ref{Lemma:3.1}, we 
have that $x > \frac{N}{2}$ so that we may conclude that $c > \frac{1}{2}$ 
(and $C > \frac{3}{2}$).\\ 
\emph{Fact:} The number of $\tau(u_j)$:s such that $\tau(u_j)\le N$ 
is equal to the 
number of $l_i:$s such that $l_i\le N$. By (\ref{cN}) this number is 
\begin{align}\label{xi}
(1-c)N.
\end{align}
 Then, for some $j\le x$, 
we must have $\delta_{j}\ge (1-c)N$. We may ask, where is this $j$? 

Define $$\rho = \rho(c,N):=\{i\mid u_i\le cN\}.$$\\

\noindent \emph{Case 1:} $j\in \rho$: Then $$C\gg \frac{\tau(u_j)}{u_j} = 
1+\frac{\delta_j}{u_j}\geq 1+\frac{1-c}{c} $$ which
is equivalent to 
$$C^2> C+1,$$ which holds if and only if $C>\Phi.$\\

\noindent \emph{Case 2:} $j\not \in \rho$:
By applying the same argument as in (\ref{cN}), we get $$\max \rho\ll c^2N.$$ 
Since $u_j>cN$, we get $\tau(u_j)= u_j+\delta_j>N$. This is equivalent to: 
For all $i$ such that $\tau(u_i)\le N$, $i\in \rho$. But then, by (\ref{xi}), 
$(1-c)N\le \max \rho$ and so, again, $C > \Phi$.
 \hfill $\Box$\\

For the rest of this section, define the sequences $a$ and $b$ as $a(1,2)$ 
and $b(1,2)$ respectively. That is, we study the solution of $(1,2)\G$.
\begin{Prop}\label{Prop:3.4} Put $R := \{b_i/a_i\mid i\in \M\}$. 
\begin{itemize}
\item[{\it a)}] The set $(\Phi, \infty ]\cap R$ is infinite.
\item[{\it b)}] Fix two constants $C \le 2\le D$ with $\beta := D - C  < 1/2$. 
 Then $([1,C)\cup (D,\infty))\cap R$ is infinite.
\item[{\it c)}] The set $[1,2]\cap R$ is infinite.
\end{itemize}
\end{Prop}

\noindent{\bf Proof.} 


\noindent {\it Item a)} Clearly $\pi_{1,2}$ satisfies the conditions of $\tau$ 
in Lemma \ref{Lemma:3.3}. The result holds since $\Phi$ is irrational 
so that for all $i$, $b_i/a_i\ne \Phi$.\\

\noindent {\it Item b)} Put $ S := \R \setminus [C, D]$. 
Suppose on the contrary that all but finitely 
many points from $R$ lie in $[C,D]$. Put 
$r := \#\{i \mid b_i/a_i\in S \}$. 
Clearly, if $b_i/a_i\in [C,D]$ with $a_i < N$, then 
$2(N-a_i) + b_i\in I(N) := [CN,DN]$. 
Denote the number of pairs $(a_i,b_i)$ with $a_i < N$ such that 
$b_i/a_i\in [C,D]$ with $J(N)$. Then,
by Lemma \ref{Lemma:3.2}, for all $N$, 
 $$J(N) > \frac{N}{2}-r.$$ 
For all $\epsilon > 0$, for all 
sufficiently large $N = N_\epsilon$, we have that 
$$\frac{J(N)-1}{N}>\frac{1}{2} - \frac{r+1}{N} > \frac{1}{2} - \epsilon.$$ 
In particular we may take $\epsilon:=1/2-D+C>0$ and choose $N$ as $N'$, 
a fixed integer strictly greater than $\frac{2(r+1)}{1-2(D-C)}$.
The number of integer points in $I(N')$ 
is $$T(N'):=\lceil DN'\rceil  - \lceil CN'\rceil.$$ Since, by 
the definition of $\epsilon $, we have 
$$\frac {\lceil DN'\rceil  - 
\lceil CN'\rceil}{N'}< 1/2-\epsilon +\frac{1}{N'}.$$ 
we get $$T(N')<J(N').$$ Then, by the Pigeonhole principle, for some $x\in I$, 
there exist a pair $i < j < N'$ such that 
$$2(N' - a_i) + b_i = x = 2(N' - a_j) + b_j.$$
But then $2(a_j-a_i)=b_j-b_i$ so that by the definition of $(1,2)\G$ 
there is a move $(a_j,b_j)\rightarrow (a_i,b_i)$, which, 
by Theorem \ref{Thm:2.1}, is impossible.\\ 

\noindent {\it Item c)} We begin by proving two claims. 
Fix an $N\in \N$ such that $b_N \ge 2a_N$. 
(Clearly there is such an $N$. Take for example $N=0$.) Then: \\

\noindent{\it Claim 1.} If there exists a least $k\in \N$ 
such that $b_{N+k} > 2a_{N+k}$, it follows that 
$b_{N+k} - 2a_{N+k} = b_N - 2a_N + 1$.
(See Table A1 for the case $N=0$ and $k=1$ and also 
the ``$\gamma$-row'' which gives 
an initial sequence of pairs ``$(N_n,k_n)$'' as follows: 
$(0,1)(1,1)(2,5)(7,1)(8,2)(10,4)(14,2)(16,k_8).)$\\

\noindent{\it Proof of Claim 1.} 
Suppose that $N>1$ is chosen smallest possible such that, contrary 
to the assumption, there is a least $k > 0$ such that 
$$\frac{b_{N+k}}{a_{N+k}} > 2,$$ and 
$b_{N+k} - b_N \ne 2(a_{N+k} - a_N) + 1.$ Then, by the minimality of $N$, 
we must have 
\begin{align}
b_{N+k} - 2a_{N+k} > b_N - 2a_N + 1.\label{Nk} 
\end{align}
But then, by the 
greedy choice of $b_{N+k}$, there must be a $j < N$ such that 
$$b_{N+k}-1-b_j = \gamma (a_{N+k}-a_j),$$ where $\gamma = 0,\frac{1}{2},1$ 
or $2$. Put $$y := b_N - 2(a_N - a_j).$$ Altogether, by (\ref{Nk}), we get 
$$b_j - y = (2-\gamma)(a_{N+k} - a_j) + 1 > 0.$$ But, by minimality of $N$, 
$b_j$ must be strictly less than $y$, a contradiction. In conclusion, 
the claim holds.\\

\noindent{\it Claim 2.} Suppose that $[1,2]\cap R$ is finite. 
Then, there is an $r\in \N$ such that $N\ge r$ implies
\begin{align}
b_{N+1} - b_N = 3 \text{ and } a_{N+1} - a_N = 1\label{31} 
\intertext{or} 
b_{N+1} - b_N = 5 \text{ and } a_{N+1} - a_N = 2.\label{52}
\end{align}

\noindent{\it Proof of Claim 2.} Since we assume that 
$[1,2]\cap R$ is finite, for some $s\in \N$, for all $j\ge s$ 
we have that $b_j / a_j > 2$. By Claim 1, since $(a_i)$ is increasing, 
this clearly implies $b_{j+1} \ge b_j + 3$. 
By definition of $a_{N+1}$, if $N$ is such that $a_N\ge b_r$, this gives
\begin{align}
a_{N+1} - a_N \le 2.\label{aj}
\end{align}
Plugging this into the result of Claim 1 
we get either $(\ref{31})$ or $(\ref{52})$. We are done.\\

The remainder of the proof consists of a geometric 
argument contradicting the greedy definition of $b$. We show 
(implicitly)  that there would be an $N$-position to much if c) 
fails to hold.

Notice that Claim 2 implies that both $a$ and $(b_i)_{r\le i}$ are
increasing. By complementarity of $a$ and $b$ it follows: $(\star )$ There 
are infinitely many $r$:s such that $(\ref{52})$ holds. 

Let $r < u < v$ be such that $(\ref{52})$ holds for both $b_u < b_v$. 
Define four lines accordingly:
\begin{align*}
l_u(x) = x+b_u,\\
l_{u+1}(x) = x+b_u+3,\\
l_v(x) = b_v,\\
l_{v+1}(x) = b_v+5,
\end{align*}
These fours lines will intersect at the integer coordinates, 
$((\alpha_i , \beta_i))_{i\in \{1,2,3,4\}} 
= ((b_v-b_u-3, b_v),(b_v-b_u, b_v),(b_v-b_u+2, b_v + 5),(b_v-b_u+5, b_v +5))$, 
defining the corners 
of a parallelogram. Denote the set of integer points \emph{strictly} 
inside this parallelogram by $\mathcal{K}$. Then, by 
inspection $$\#\mathcal{K} = 8$$ and, by $(\star )$, 
we may assume that we have chosen $v$ sufficiently large so that, for all 
$(x,y)\in \mathcal{K}$, 
\begin{align}
 1 < y/x < 2.\label{yx} 
\end{align}

Denote by $\mathcal{L}$ another set of lines satisfying 
the following conditions. A line $l$ belongs to $\mathcal{L}$ if and only if: 
\begin{enumerate}[(i)]
\item Its slope is either $1/2, 2$ or $\infty$. 
\item It intersect a point of the form $(a_s,b_s)$ or $(b_s,a_s)$ 
with $s\ge r$.
\item It intersects $\mathcal{K}$.
\end{enumerate}

By the definition of $\mathcal{K}$ it follows from $(\ref{52})$ 
that $i$ and $j$ may be defined such that 
each line of form (ii) and (iii) is also of the form (i). 
Again, by $(\star )$ we may assume that we have chosen $i$ and $j$ 
sufficiently large so that the first part of (ii) together with (iii) 
implies $s\ge r$. \\

\noindent{\it Claim 3:} There is an integer coordinate in the set 
$\mathcal{K}\setminus \mathcal{L}$.\\

Clearly, by the definition of $(b_i)$ and by (\ref{yx}), 
this claim contradicts the assumption that $[1,2]\cap R$ is finite. 
(In fact it would imply the existence of an $N$-position in $\mathcal{K}$ 
without a $P$-position as a follower.)\\ 

\noindent{\it Proof of Claim 3.} Let 
$$\mathcal{K}' := \{(0,0),(1,0),(1,1),(2,1),(2,2),(3,2),(3,3),(4,3)\}.$$ 
Then $\mathcal{K}'$ is simply 
a linear translation of $\mathcal{K}$. (Namely, given $(x,y)\in \mathcal{K}$, 
$T(x, y) = x - (b_j - b_i - 1), y - (b_j + 1)\in \mathcal{K}'$.) 

Let $\alpha \in \R$. Clearly, the two lines $x + \alpha$ and 
$x + 3 + \alpha$ can together cover at most three points in $\mathcal{K}'$, 
namely choose $\alpha = 0$ or $1$. The two lines $2x - \alpha$ and 
$2x - 5 - \alpha$ can cover at most two points in $\mathcal{K}'$, namely we may 
choose $\alpha = 0, 2$ or $3$. (In fact, for the two latter cases it is 
only the former line that contributes.) On the other hand, the 
two lines $x/2 + \alpha$ and $x/2 + 5/2 + \alpha$ can cover at most two points 
in $\mathcal{K}'$, namely, if we choose $\alpha = 0, 1/2$ or $1$. (In fact, 
 for these $\alpha$, it is only the former line that contributes.) 

Fix any set of the above six lines, depending only on the choices of $\alpha$ 
for the respective cases, and denote this set by $\mathcal{L}'$. Then, as 
we have seen,  $\#(\mathcal{L}'\cap \mathcal{K}')\le 7$.  
But, by Claim 2, an instance of $\mathcal{L}\cap \mathcal {K}$ 
is simply a linear translation of some set $\mathcal{L}'\cap \mathcal{K}'$. 
The claim follows and so does the proposition. $\hfill \Box$

\begin{Rem}
Obviously, Proposition \ref{Prop:3.4} a) may be adapted without 
any changes for general $p$ and $q$. Also, item b) may easily 
be generalized. On the other hand, we did not find any immediate way 
to generalize c) in its present form. 
\end{Rem}

\begin{Thm}\label{Thm:3.1} Let $a = a(1,2)$ and $b = b(1,2)$. Then the limit 
$$\lim_{n\in \N}\frac{b_n}{a_n}$$ does not exist.
\end{Thm}
\noindent{\bf Proof.} 
Suppose on the contrary that $\alpha := \lim_{n\in \N}\frac{b_n}{a_n}$ exists.
Then either
\begin{enumerate}[(i)]
\item $\alpha \in [1,\Phi)$, 
\item $\alpha \in [\Phi,2]$, or
\item $\alpha \in (2,\infty]$.
\end{enumerate}
By Proposition \ref{Prop:3.4} a), (i) is impossible. On the other 
hand (ii) is contradicted by Proposition \ref{Prop:3.4} b) with, say, $C=\Phi$ 
and $D=2$. For the last case Proposition \ref{Prop:3.4} c) 
gives a contradiction. $\hfill \Box$\\

\begin{Cor}
Define $a$ and $b$ as in Theorem \ref{Thm:3.1}. Then: 
\begin{enumerate}[(i)]
\item For all $n\in \N$ there exist $ i,j\ge n$ such that $$\beta = 
\left|\, \frac{b_i}{a_i}-\frac{b_j}{a_j}\,\right| \ge \Phi -\frac{3}{2}.$$
\item If $((a_i,b_i))$ splits and the conditions for $(l_i)$ and $(u_i)$ 
in Conjecture \ref{Conj:3} are satisfied, then $1<\eta <\Phi < \gamma \le 3$.
\end{enumerate}
\end{Cor}
\noindent{\bf Proof.} Item (i) is a consequence of Proposition \ref{Prop:3.4}. 
If, for infinitely many $i$:s,  $b_i/a_i\le 3/2$, then by b), 
there has to be infinitely many $j$:s such that $b_j/a_j\ge \Phi$. On 
the other hand, if there are not infinitely many $i$:s of the first form, 
then, by c), there has to be infinitely many $i$:s such that $b_i/a_i\ge 2$. 
Then, again, by c) we may choose $C=2$ and $D=5/2$, 
which implies $\beta \ge 1/2$. 
For item (ii), by (i) and Proposition \ref{Prop:3.4}, it only remains 
to verify that $\gamma \le 3$. But this follows by the greedy choice 
of $\pi_{1,2}$, since the worst case is if, for all but finitely 
many $i$, $b_i/a_i > 2$. But then, Lemma \ref{Lemma:3.2} gives the result. 
\hfill $\Box$\\

In \cite{Lar2} a restriction, called Maharaja Nim, of the game $(1,2)$GDWN 
is studied. Here, all options on the two $(1,2)$ diagonals, 
except $(1,2)$ and $(2,1)$, are forbidden. 
In contrast to the main result of this paper, for Maharaja Nim 
it is proved that the $P$-positions lie on the 
same 'beams' as in Wythoff Nim, 
however the behaviour along these beams turns out to be fairly 'chaotic'.
In this context it is interesting to observe that the only option on 
the $(1,2)$ diagonal which is a splitting pair is $(1,2)$ itself. This 
stands in bright contrast to the main result, Proposition \ref{Prop:4.3}, 
of the final section in this paper. 
Namely, for $(1,2)$GDWN the non-splitting pairs on 
the diagonal $(1,2)$ contributes 
significantly in destroying the asymtotics of the $P$-positions 
of Maharaja Nim.

\section{More on splitting pairs}
Let $G$, $H$ be impartial games. Then, if 
$\mathcal{P}(G) = \mathcal{P}(H)$, 
we say that $G$ is \emph{equivalent} to $H$. 
The main result of this section is a partial resolution of Conjecture 
\ref{Conj:1}. 
(See also for example \cite{BFG, DFNR, FHL} for related results.)

\begin{Prop}\label{Prop:4.3}
Suppose that $(p, q)$ is a non-splitting pair with $1 < \frac{q}{p} < \Phi$. 
Then $(p,q)\G $ is equivalent to Wythoff Nim.
\end{Prop}

 Before proving this proposition we need to develop some facts 
from combinatorics on Sturmian words. 
We will make use of some 
terminology and a lemma from \cite[Section 1 \& 2]{Lot01}. 

Let us define two infinite words $s$ and $s'$. For all $n\in \M$, 
the $n$:th letter is 
$$s(n) := \lfloor (n+1)\Phi \rfloor - \lfloor (n)\Phi \rfloor$$ and 
$$s'(n) := \lceil (n+1)\Phi \rceil - \lceil (n)\Phi \rceil,$$
respectively. Then $s$ is the \emph{lower mechanical} word with 
\emph{slope} $\Phi$ and \emph{intercept} $0$ and $s'$ 
is the \emph{upper} ditto. Whenever we want to emphasize that $\Phi$ is 
irrational we say that $s$ (or $s'$) is \emph{irrational mechanical}. 
Thus, the \emph{characteristic} word belonging 
to $s$ and $s'$ is $c = s(1)s(2)s(3)\ldots$. Namely, we have  
$s(0) = 1$, $s'(0) = 2$ and otherwise, for all $n>0$, 
$s(n) = s'(n) = 1$ or $s(n) = s'(n) = 2$. In fact, we have 
$$s=12122121221\ldots $$ and $$s'=22122121221\ldots .$$
Denote with $ l(x)$ the number of letters in $x$ 
and with $h(x)$ the number of $1$:s in $x$. 
Let $\alpha $ and $\beta $ be two \emph{factors} of a word $w$. Then $w$ 
is \emph{balanced} if 
$l(\alpha) = l(\beta)$ implies $\mid h(\alpha)- h(\beta)\mid \ \le 1$. 
By \cite[Section 2]{Lot01}, both $s$ and $s'$ are 
balanced (aperiodic) words. We will also need the following result.  

\begin{Lemma}[\cite{Lot01}]\label{Lemma:4.1} 
Suppose two irrational mechanical words have the same slope. 
Then their  respective set of factors are identical.
\end{Lemma}
We also use the following notation.
Suppose $x = x_1x_2\ldots x_n$ is a factor of a mechanical word on $n$ letters. 
Then we define the \emph{sum} of $x$ as $\sum x := x_1+x_2+\ldots +x_n$. 
For example the sum of $2121$ equals $6$.
We let $\xi_n(s)$ denote the unique \emph{prefix} of an infinite 
word $s$ on $n\in \M$ letters.
\begin{Lemma}\label{Lemma:4.2}
Let $x$ be any factor of $s$ (or $s'$). 
Then $$\sum x = \sum \xi_{l(x)}(s) \text{ or } \sum x=\sum \xi_{l(x)}(s').$$ 
\end{Lemma}
\noindent {\bf Proof.} If two factors of $s$ have the same length and the 
same height, then, since the number of 2:s in the respective factors must 
be the same, their sums are identical. Therefore, 
if $h(x) = h(\xi_{l(x)}(s))$, this implies $\sum x = \sum\xi_{l(x)}(s)$. 

Assume on the contrary that $h(x)\ne h(\xi_{l(x)}(s))$. 
On the one hand, for all $n$, $ h(\xi_{n}(s)) = h(\xi_{n}(s')) + 1$. 
On the other hand, 
the balanced condition implies: Given an $l'\in \N$ if $x$ is a factor of $s$ 
such that $l(x)=l'$ then $h(x)$ is one of two fixed values.  
It follows that $h(x) = h(\xi_{l(x)}(s'))$. But then, by the initial 
observation, we are done. 
$\hfill \Box $\\

The following proposition assures that for each splitting pair, a 'split 
is initiated'. 

\begin{Prop}\label{Prop:4.4}
Let $p,q\in N$. Then $(p,q)$ is a splitting pair if and only if 
there exists a pair $m, n\in \M$ 
with $m < n$ such that 
$$(p,q)=(a_n-a_m,b_n-b_m)=(\lfloor n\Phi\rfloor - \lfloor m\Phi\rfloor, 
\lfloor n\Phi^2\rfloor - \lfloor m\Phi^2\rfloor).$$
\end{Prop}

\noindent {\bf Proof.} Suppose that $(p,q)$ is a splitting pair. If 
$(p,q)=(a_n,b_n)$, for some $n\in \N$, we may take $m=0$. 
If $(p,q)=(a_r+1,b_r+1)$, for some $r\in \N$, then, 
since $s$ and $s'$ are mechanical with the same slope, 
by Lemma \ref{Lemma:4.1} $\rho = \xi_r(s')$ 
is a factor of $s$. But then,
$a_r+1 = \sum \rho = \sum \xi_n(s) - \sum \xi_m(s)$ for some $n-m=r$. 
This gives $(p,q)=(a_r+1,a_r+1+r)=(a_n-a_m, a_n-a_m+n-m) = (a_n-a_m, b_n-b_m)$.
For the other direction, let $(p,q)$ and $m<n$ be as in the proposition. 
Then, define $x:=\xi_n(s)-\xi_m(s)$, and so, by $l(x)=n-m$ and 
Lemma \ref{Lemma:4.2},  
we may take $p=\sum\xi_{l(x)}(s)=a_{l(x)}$ or $p=\sum\xi_{l(x)}(s')=a_{l(x)}+1$. 
In either case, the assumption gives $q=p+l(x)$, and so 
$(p,q)$ is a splitting pair. $\hfill \Box $\\

\noindent {\bf Proof of Proposition \ref{Prop:4.3}.} We need to show that, 
for all $i$, 
$(a^{1,1}_i,b^{1,1}_i) = (a^{p,q}_i, b^{p,q}_i)$. Let $i\in \M$ and denote 
with $(a_i,b_i) = (a^{1,1}_i,b^{1,1}_i)$. By Proposition \ref{Prop:4.4}, 
there is a pair $j < i$ such that $(a_i,b_i)\rightarrow (a_j,b_j)$ is 
a legal move of $(p,q)\G$ if and only if $(p,q)$ is a splitting pair. 

Suppose that there is a 
pair $j < i$ such that $(a_i,b_i)\rightarrow (b_j,a_j)$ is a legal move. 
Then $(a_i-b_j, b_i-a_j) = t(p,q)$ for some $t\in  \N$. Hence $$a_i=tp+b_j$$
 and $$b_i=tq+a_j.$$ Then $\Phi < b_i/a_i = (tq+a_j)/(tp+b_j) < q/p$ 
since $0\le a_j\le b_j$ for all $j$. 
 $ \hfill \Box $\\

\appendix
Here we present various tables and figures representing $P$-positions 
of $Q\G $  for different $Q$. We have also plotted the ratios $a_n/b_n$ 
on intervals $[0,n]$, for $n$ up to $50000$. (The code is written in C, 
in fact the original code stems from Jonas Knape's and my Master's Thesis.) 
The purpose of this appendix 
is to support our conjectures and stimulate further questions and research 
on generalized Wythoff games.
\vspace{1 cm}
\begin{table}[ht!]
\begin{small}
\begin{center}
\setlength{\tabcolsep}{3pt}{
  \begin{tabular}
{| l || c|c|c|c| c| c| c| c| c| c| c| c| c| c| c| c| c| c| c|}
    \hline
    $b_n$\T   &0&\bf{3}&\bf{6}&5&10&14&17&\bf{25}&\bf{28}&18&\bf{35}&23&31&29&\bf{48}&32&\bf{55}&37&40\\
    $a_n$\T   &0&1&2&4&7 &8 &9 &11&12&13&15&16&19&20&21&22&24&26&27\\\hline
 $\delta_n$\T &0&2&4&1&3 &6 &8 &14&16&5 &20&7 &12&9 &27&10&31&11&13\\
$\gamma_n$\T  &0&1&2&-3&-4&-2&-1&3&4&-8&5&-9&-7&-11&6&-12&7&-15&-14\\
$\eta_n$\T  &0&5&10&6&13&20&25&39&44&23&55&30&43&38&75&42&86&48&53\\\hline
    $n$ \T    &0&1&2&3&4 &5 &6 &7 & 8&9 &10&11&12&13&14&15&16&17&18\\
    \hline
  \end{tabular}}
\end{center}
\end{small}
\caption{Here $\{a_n, b_n\}$ represents a $P$-position 
of $(1, 2)\G$ for $0\le n\le 18$. 
Further, $\delta_n = b_n - a_n$, $\gamma_n = b_n - 2a_n, \eta_n = 2b_n - a_n$.}
\end{table}

\begin{table}[ht!]
\begin{center}
  \begin{tabular}
{| l || c|c|c| c| c| c| c| c| c| c| c| c| c|}
    \hline
    $b_n$\T   &0&2&6&8&7 &16&18 &20&17&24&26&21&34\\  
    $a_n$\T   &0&1&3&4&5 &9 &10 &11&12&13&14&15&19\\\hline  
 $\delta_n$\T &0&1&3&4&2 &7 &8  &9 &5 &11&12&6 &15\\  
    $n$ \T    &0&1&2&3&4 &5 &6  &7 &8 &9 &10&11&12\\  
    \hline
  \end{tabular}
\end{center}\caption{The first $P$-positions of $(2,3)\G$ and 
$\delta_n = b_n - a_n$.}
\end{table}

\begin{table}[ht!]
\vspace {1 cm}
\begin{center}
  \begin{tabular}
{| l || c|c|c| c| c| c| c| c| c| c| c| c|}
    \hline
    $b_n$\T   &0&2&5&7&10 &17&14&19&18&20&27&33\\  
    $a_n$\T   &0&1&3&4&6  &8 &9 &11&12&13&16&21\\\hline  
 $\delta_n$\T &0&1&2&3&4 &9 &5  &8 &6 &7&11&12 \\  
    $n$ \T    &0&1&2&3&4 &5 &6  &7 &8 &9 &10&11\\  
    \hline
  \end{tabular}
\end{center}\caption{Here $\{a_n,b_n\}$ represents a 
$P$-position of $(2,4)\G$. 
Notice that $(8,13) \ominus (2,1) = 3\times(2,4),$ so that $(8,13)$ is the 
first Wythoff-pair that short-circuits $(2,4)\G$. So, a 
`split is initiated', but our computations suggest that  
 the quotient $b_n/a_n$ converges to $\Phi$ (see Conjecture \ref{Conj:1}.
}\label{A3}
\end{table}

\begin{table}[ht!]
\vspace {1 cm}
\begin{center}
  \begin{tabular}
{| l || c|c|c| c| c| c| c| c| c| c| c| c|}
    \hline
    $b_n$\T   &0&2&5&7&10 &17&14&19&18&20&27&33\\  
    $a_n$\T   &0&1&3&4&6  &8 &9 &11&12&13&16&21\\\hline  
 $\delta_n$\T &0&1&2&3&4 &9 &5  &8 &6 &7&11&12 \\  
    $n$ \T    &0&1&2&3&4 &5 &6  &7 &8 &9 &10&11\\  
    \hline
  \end{tabular}
\end{center}\caption{The first $P$-positions of $(4,6)\G$.}
\end{table}

\begin{table}[ht!]
\vspace{1 cm}
\begin{center}
  \begin{tabular}
{| l || c|c|c| c| c| c| c| c| c| c| c| c|}
    \hline
    $b_n$\T   &0&2&5&8&9 &13 &17 &16&20&25&24&28\\  
    $a_n$\T   &0&1&3&4&6  &7 &10 &11&12&14&15&18\\\hline  
 $\delta_n$\T &0&1&2&4&3 &6 &7  &5 &8 &11&9&10 \\  
    $n$ \T    &0&1&2&3&4 &5 &6  &7 &8 &9 &10&11\\  
    \hline
  \end{tabular}
\end{center}\caption{The first $P$-positions of $(4,7)\G$.}
\end{table}
\begin{table}[ht!]
\begin{center}
\vspace{1 cm}
\begin{tabular} {c c c c c c c c c c c c } 
1, & 2 & 3, & 5 & 8, & 13 & 21, & 34 & 55, & 89 & \ldots \\ 
4, & 7 & 11, & 18 & 29, & 47 & 76, & 123 & 199, & 322 & \ldots \\ 
6, & 10 & 16, & 26 & 42,  & 68  & 110, & 178  & 288,  & 466 & \ldots \\ 
9, & 15  & 24, & 39 & 63, & 102 & 165, & 267 & 432, & 699  & \ldots \\ 
12, & 20  & 32, & 52 & 84, & 136 & 220, & 356 & 576, & 932  & \ldots \\ 
\end{tabular}
\end{center}\caption{The first five rows of the Wythoff Array. 
The Wythoff pairs are pairs of entries of the form $x,y$. For example 
9, 15 is a Wythoff pair, but 15 24 is not. On the other hand 15, 24 is 
a Dual Wythoff pair.}
\begin{center}
\vspace{5 mm}
\begin{tabular} {c c c c c c c c c c c c } 
2, & 3 & 5, & 8 & 13, & 21 & 34, & 55 & 89, &  144 &\ldots \\ 
4, & 6 & 10, & 16 & 26, & 42 & 68, & 110 & 178, & 288 & \ldots \\ 
7, & 11 & 18, & 29 & 47,  & 76  & 123, & 199  & 322,  & 521 & \ldots \\ 
9, & 14 & 23, & 37 & 60, & 97 & 157, & 254 & 311, & 565  & \ldots \\ 
12, & 19 & 31, & 50 & 81, & 131 & 212, & 343 & 555, & 878  & \ldots \\ 
15, & 24 & \vdots & \vdots & \vdots & \vdots & \vdots & \vdots 
& \vdots & \vdots  & $\ddots $\\ 
\end{tabular}
\end{center}\caption{The Dual Wythoff Array. 
The Dual Wythoff pairs are pairs of entries of the form $x,y$. For example 
9, 14 is a Dual Wythoff pair, but 14 23 is not. On the other hand 14, 23 is 
a Wythoff pair, namely the first two entries in the sixth row. This follows 
because 14 is the least number not contained in the first five rows of the 
Wythoff array and the second entry may be defined via the so-called Zeckendorff 
right-shift $Z$ of $14 = 1 + 13$ in the first row, 
namely $Z(1+13) = 2+21 = 23$.} 
\end{table}

\clearpage
\begin{figure}[ht!]
\begin{center}
	  {\includegraphics[width=0.5\textwidth]{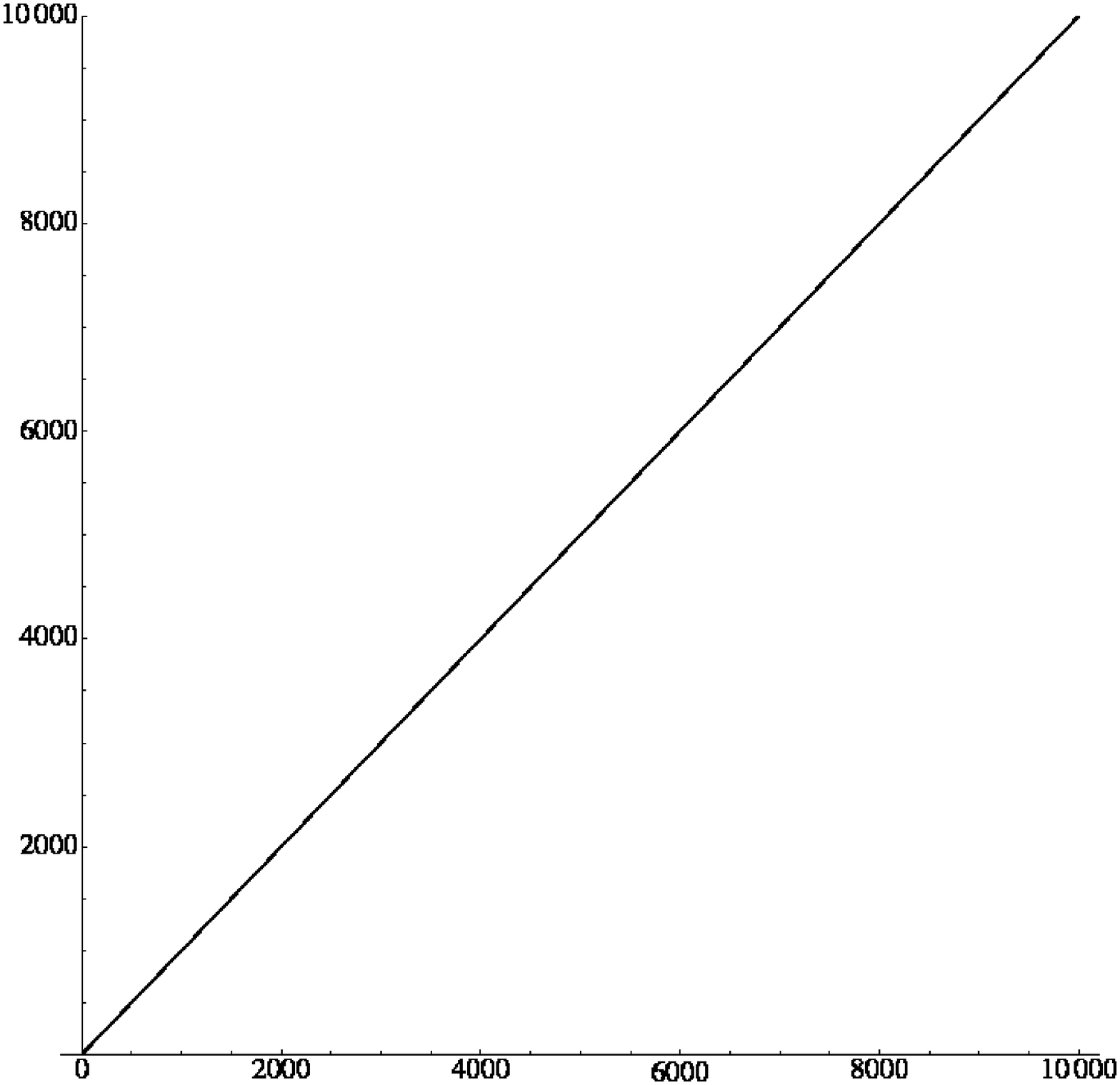}}  
\end{center}\caption{The positions $(n,\pi_Q(n))$ for $Q=\{(0,1)\}$. 
The single ``Nim beam''.}
\begin{center}
	  {\includegraphics[width=0.5\textwidth]{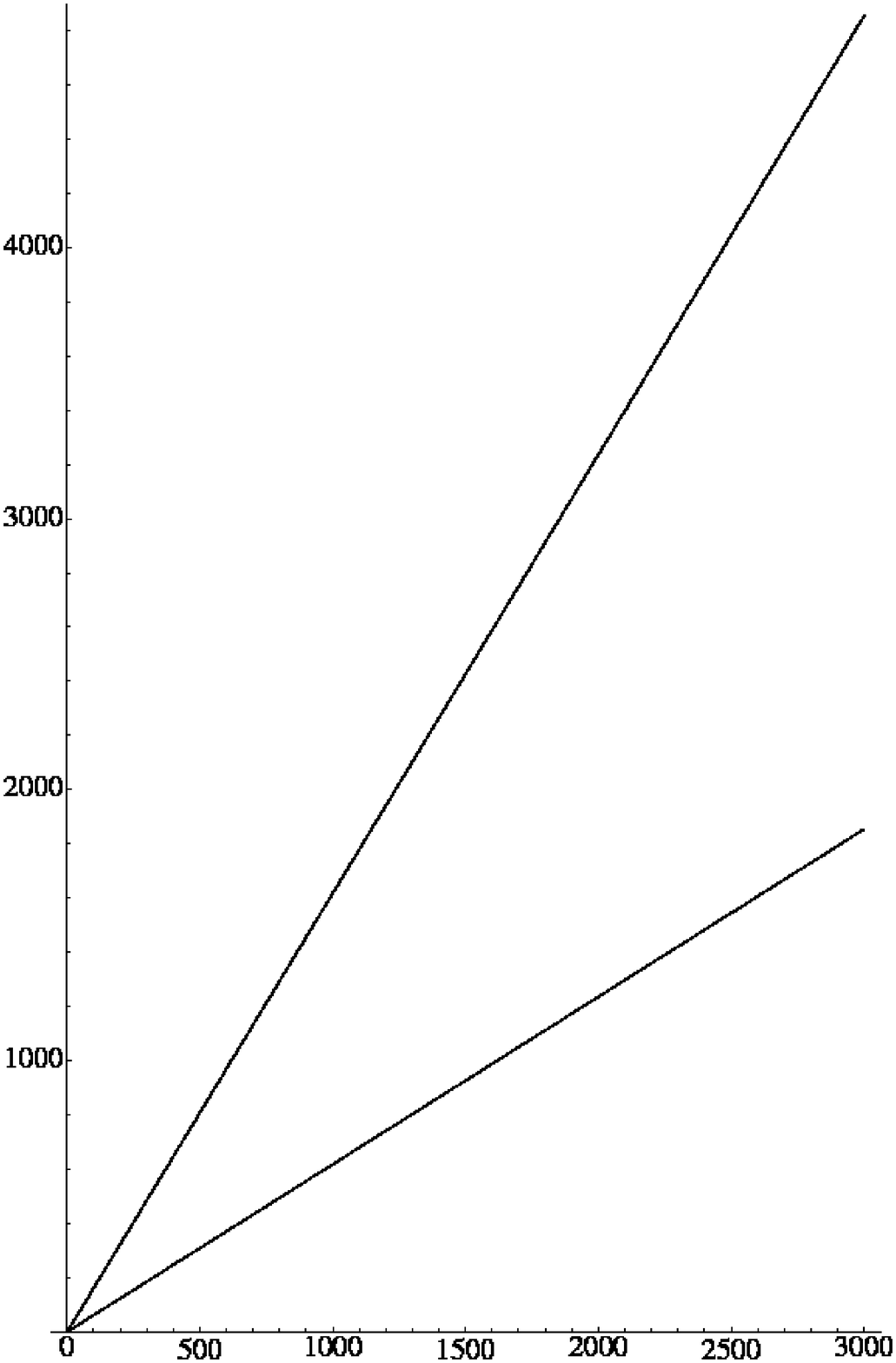}}  
\end{center}\caption{The positions $(n, \pi_Q(n))$ for $Q=\{(0,1),(1,1)\}$. 
The fundamental $1$-split of Wythoff Nim.}
\end{figure}
\clearpage
\begin{figure}[ht!]
\begin{center}
	  {\includegraphics[width=0.65\textwidth]{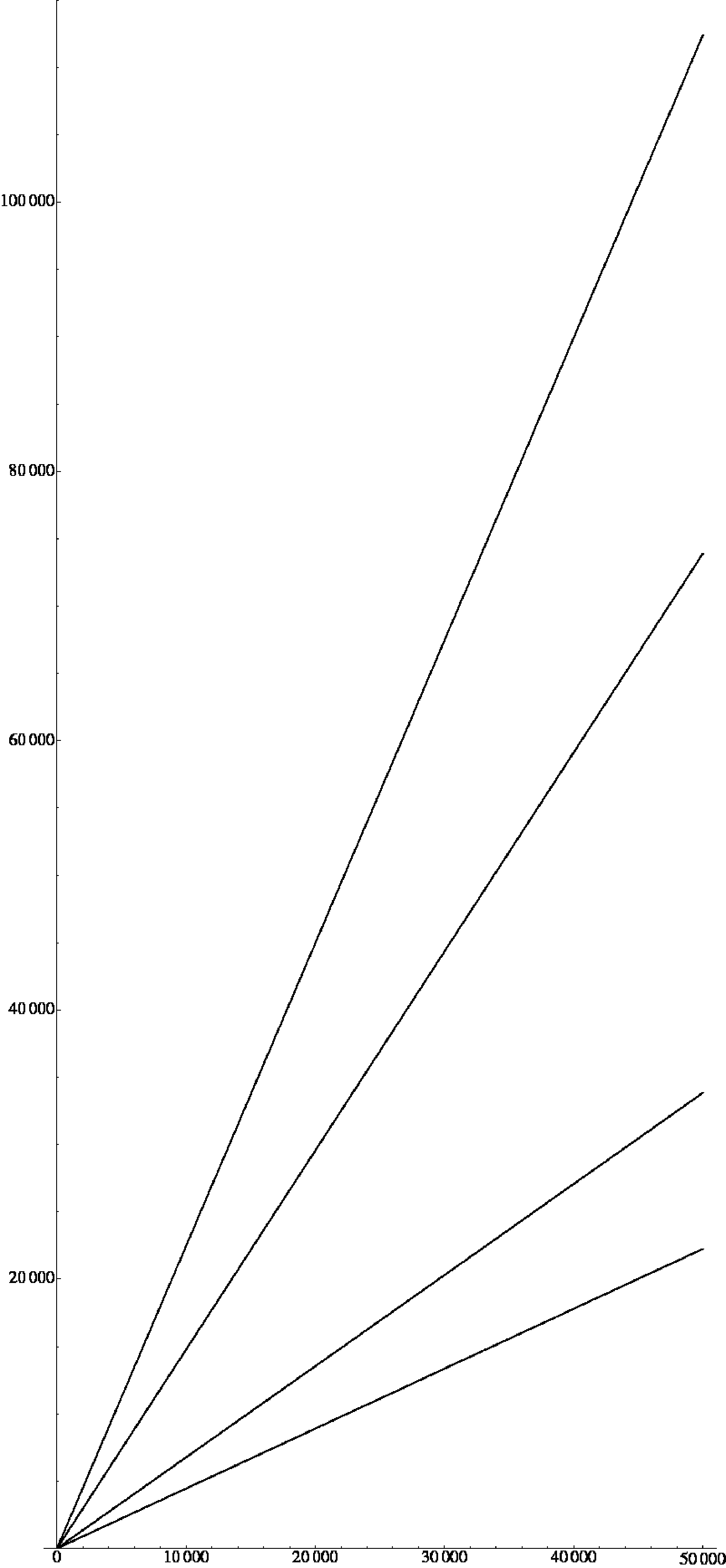}}  
\end{center}\caption{The $P$-positions $(n, \pi_Q(n))$ of $(1,2)\G$ 
and $0\le n \le 50000$.}
\end{figure}
\clearpage
\begin{figure}[ht!]
\begin{center}
	  {\includegraphics[width=0.8\textwidth]{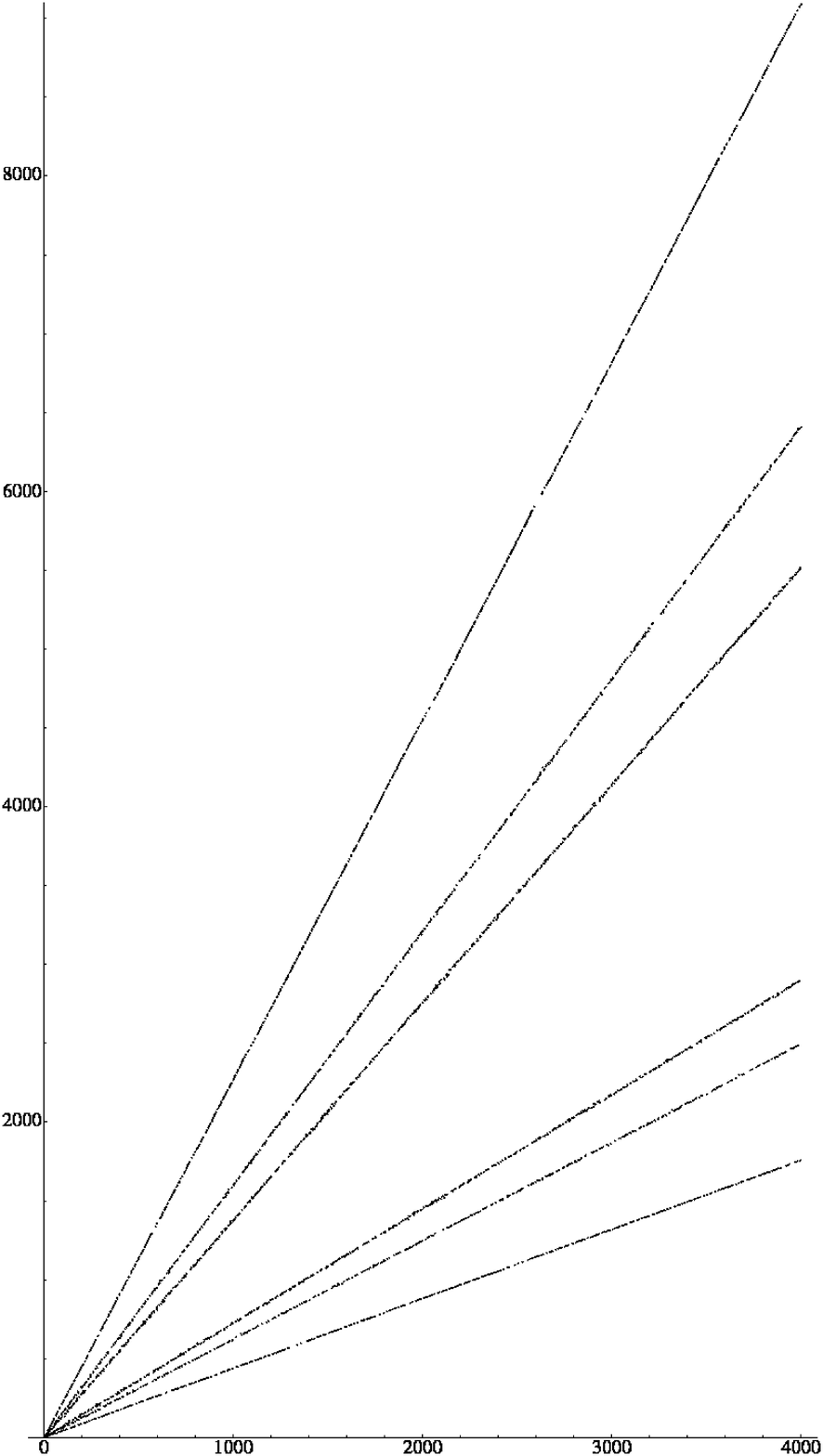}}  
\end{center}\caption{The $P$-positions $(n,\pi_Q(n))$ of $\{(1,2),(2,3)\}\G$ 
and $0\le n \le 4000$.}
\end{figure}
\clearpage
\begin{figure}[ht!]
\begin{center}
	  {\includegraphics[width=0.8\textwidth]{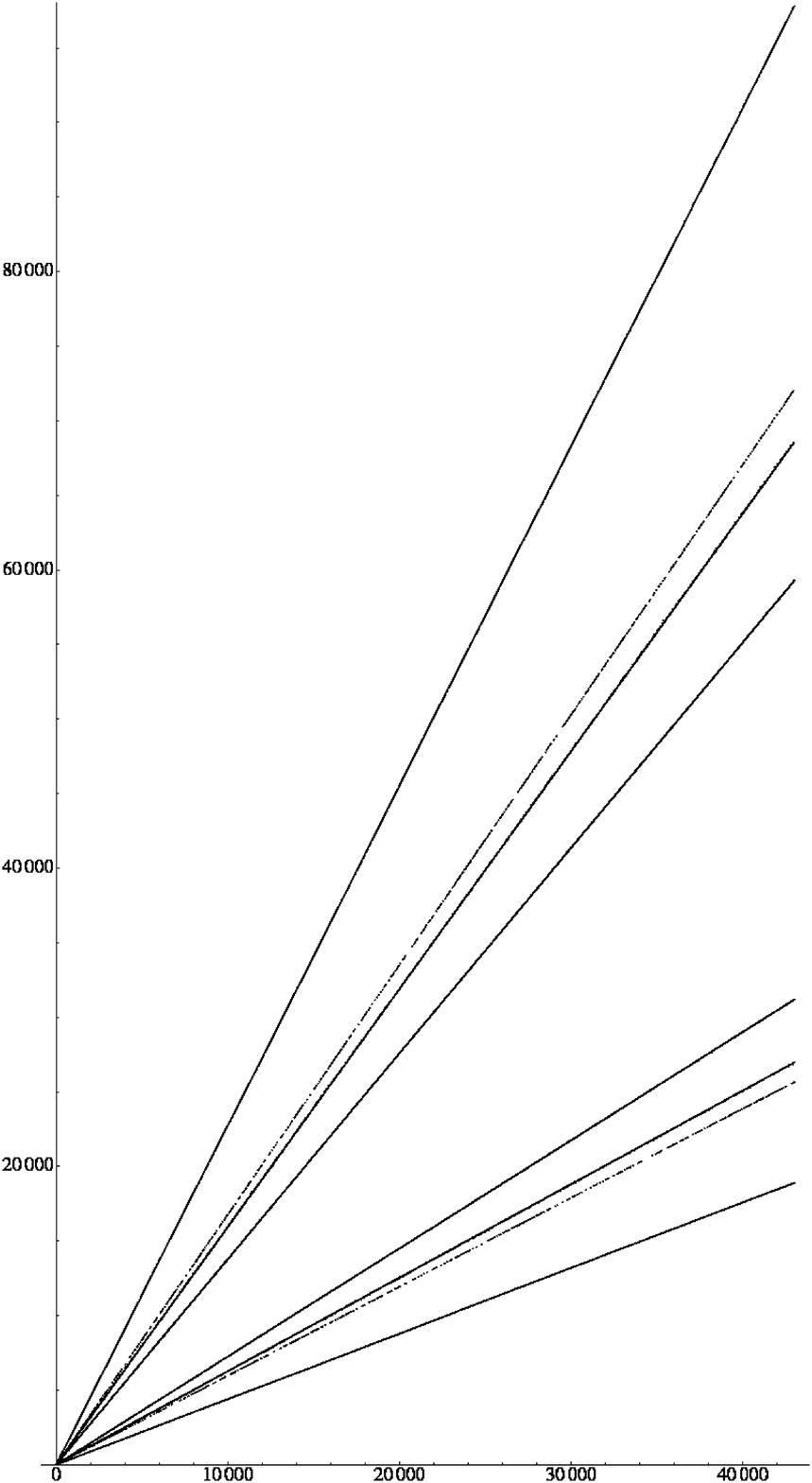}}  
\end{center}\caption{The $P$-positions  $(n,\pi_Q(n))$ 
for $\{(1,2),(2,3),(3,5)\}\G$ and $0\le n \le 40000$.}
\end{figure}
\clearpage
\begin{figure}[ht!]
\begin{center}
	  {\includegraphics[width=0.5\textwidth]{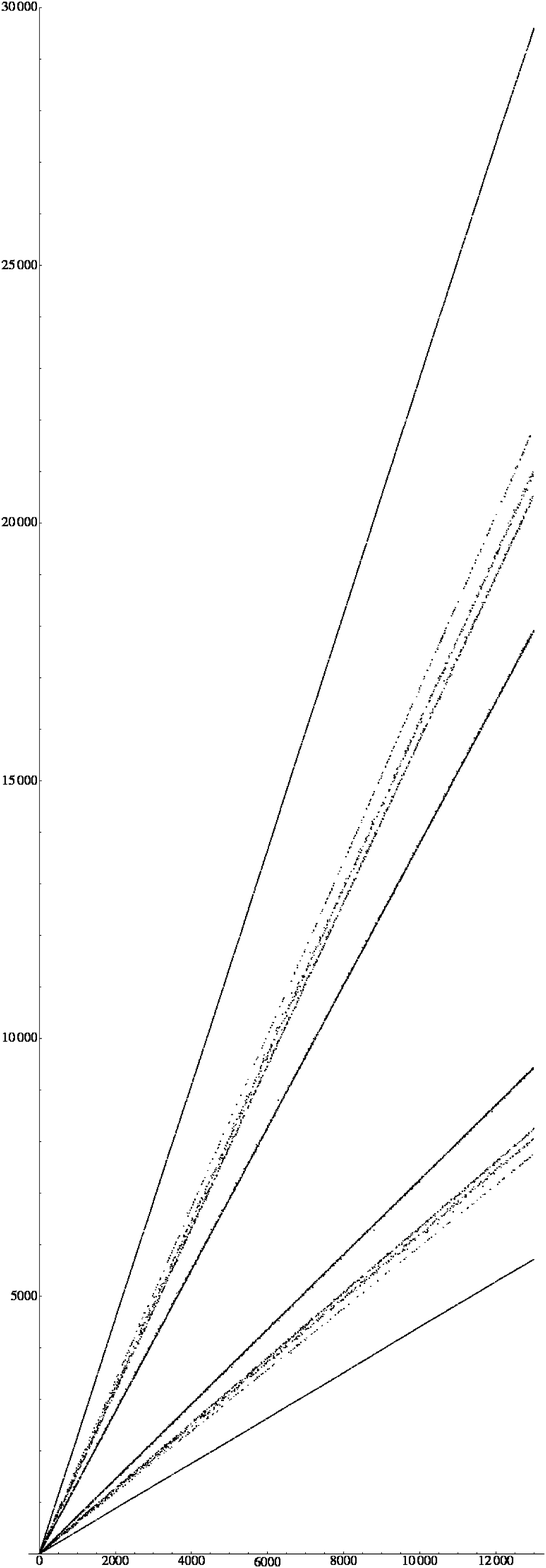}}  
\end{center}\caption{The $P$-positions $(n,\pi_Q(n))$ of
$\{(1,2),(2,3),(3,5),(5,8)\}\G$ and $0\le n \le 10000$.}
\end{figure}
\clearpage
\begin{figure}[ht!]
\begin{center}
	  {\includegraphics[width=0.7\textwidth]{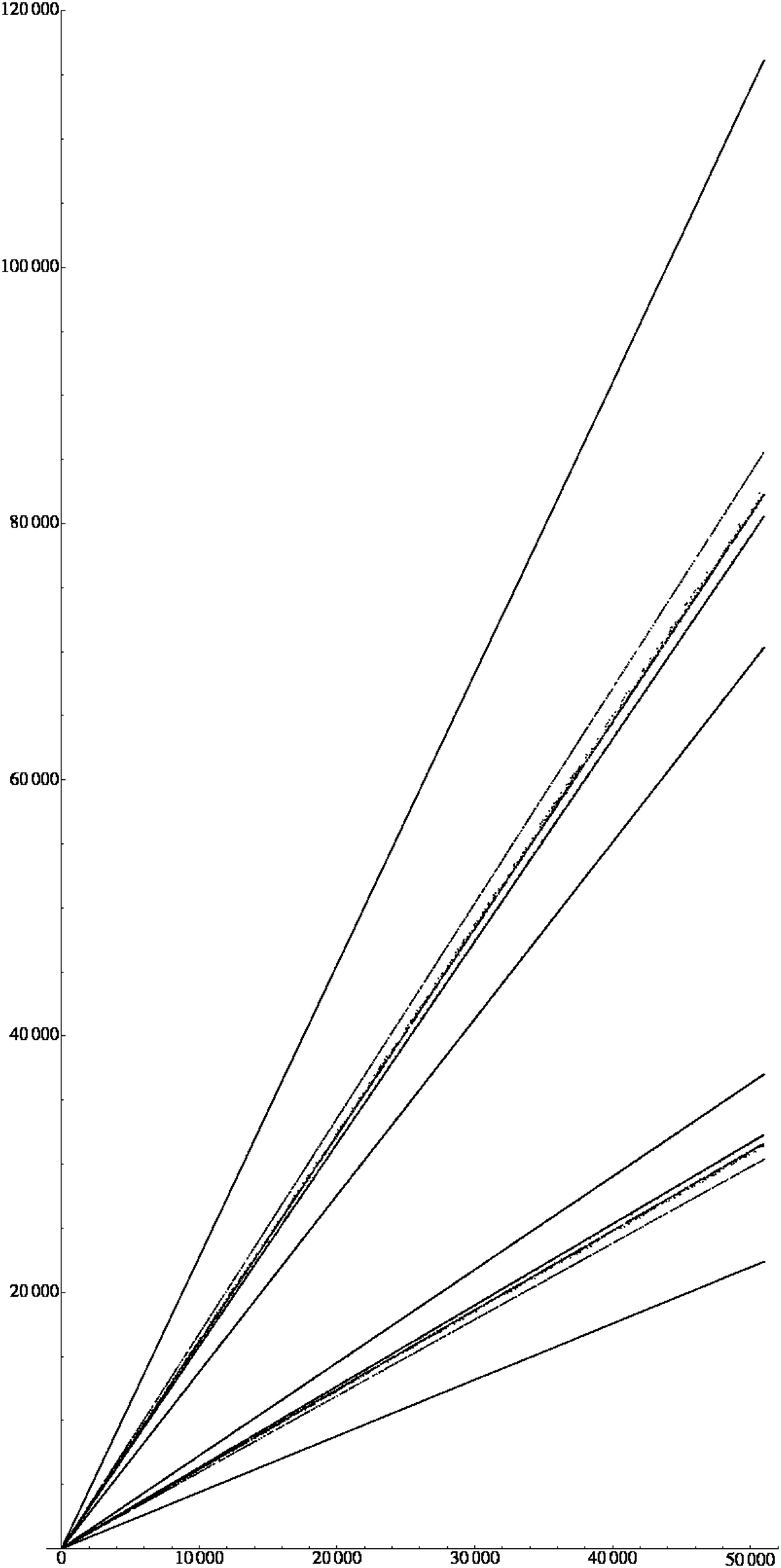}}  
\end{center}\caption{The $P$-positions $(n,\pi_Q(n))$ of
$\{(1,2),(2,3),(3,5),(5,8),(8,13)\}\G$ and $0\le n \le 51000$.}
\end{figure}
\clearpage
\begin{figure}[ht!]
\begin{center}
    {\includegraphics[width=0.7\textwidth]{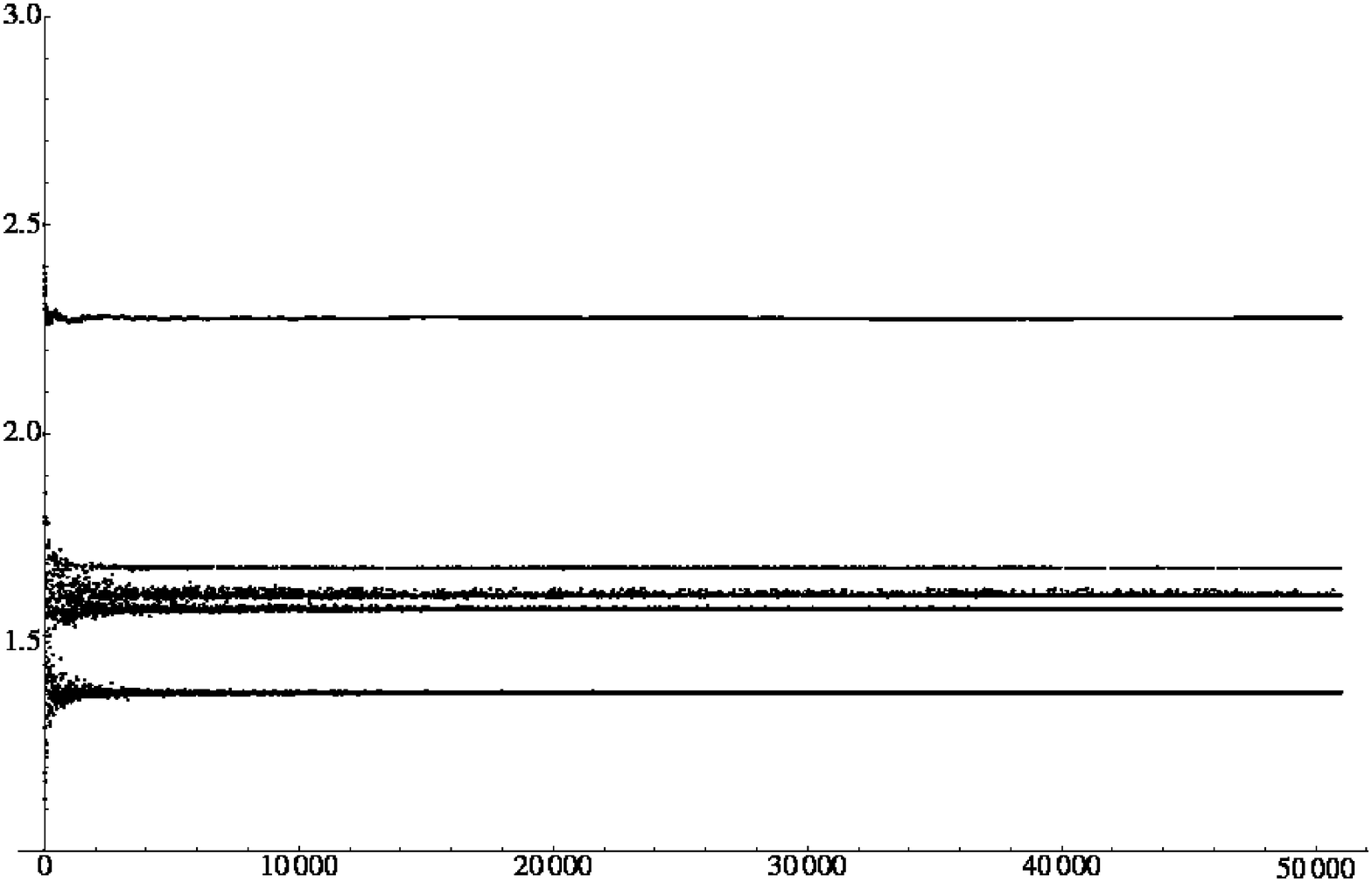}}  
\end{center}\caption{The ratio $\frac{b_n}{a_n}$ of
$\{(1,2),(2,3),(3,5),(5,8),(8,13)\}\G$ and $0\le n \le 51000$.}\label{A8}
\end{figure}
\begin{figure}[ht!]
\begin{center}
	  {\includegraphics[width=0.7\textwidth]{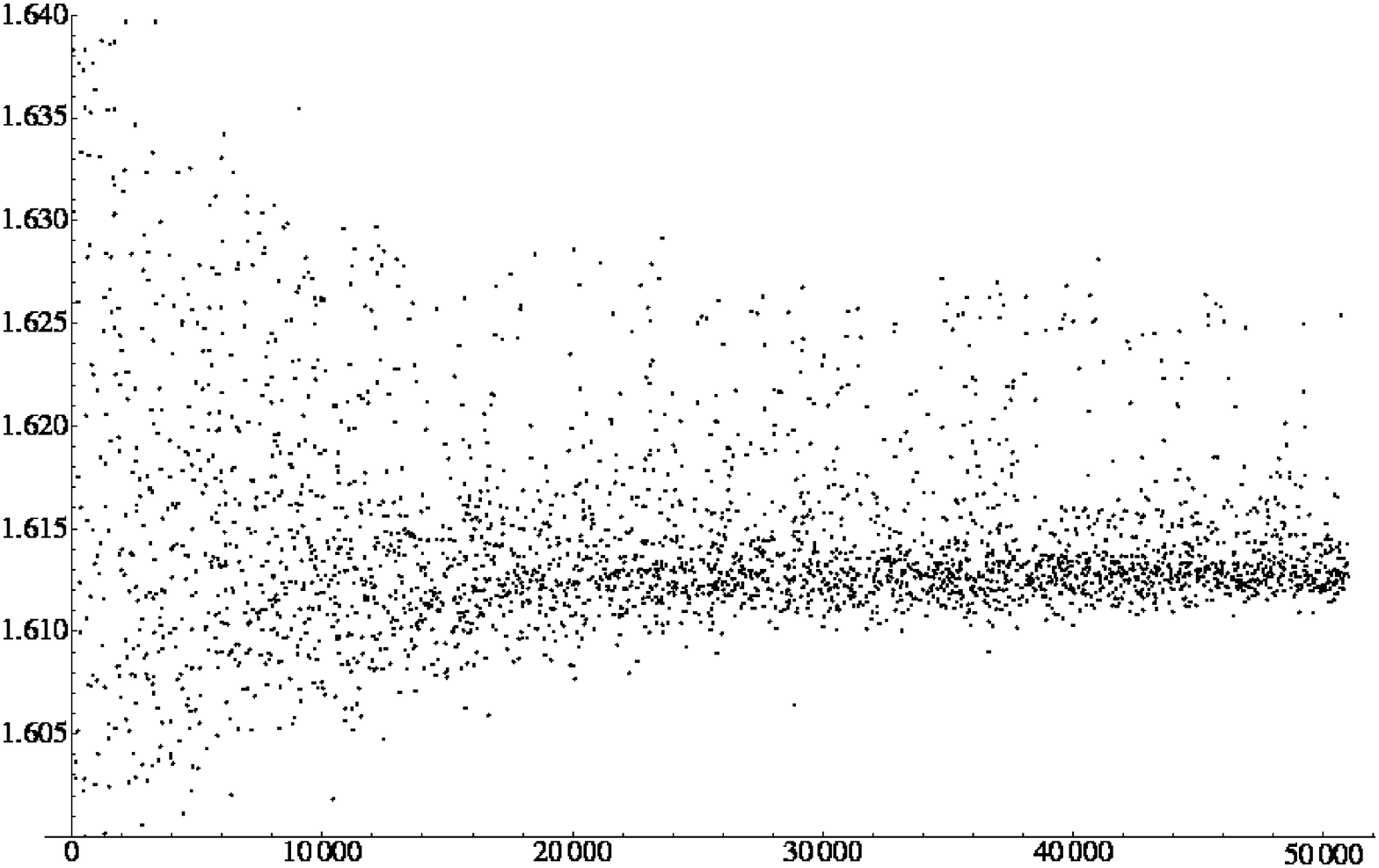}}  
\end{center}\caption{The ratio $\frac{b_{n_i}}{a_{n_i}}$ of the 'upper 
central beam' for $\{(1,2),(2,3),(3,5),(5,8),(8,13)\}\G$ and 
$0\le n \le 51000$. Does this 'perturbed beam' eventually split into two 
new beams one above $\Phi$ and the other below $\Phi$? Is there some $(p,q)$ 
that 'splits the $P$-positions' of $\{(1,2),(2,3),(3,5),(5,8)\}\G$ further 
into 6 beams above the main diagonal? Of course, we do not even know whether 
$\{(1,2),(2,3),(3,5),(5,8)\}\G$ splits into 5 beams, but indeed, 
by Figure \ref{A8} our experimental data suggests that this may hold.}
\end{figure}
\clearpage
\begin{figure}[ht!]
\begin{center}
	  {\includegraphics[width=0.7\textwidth]{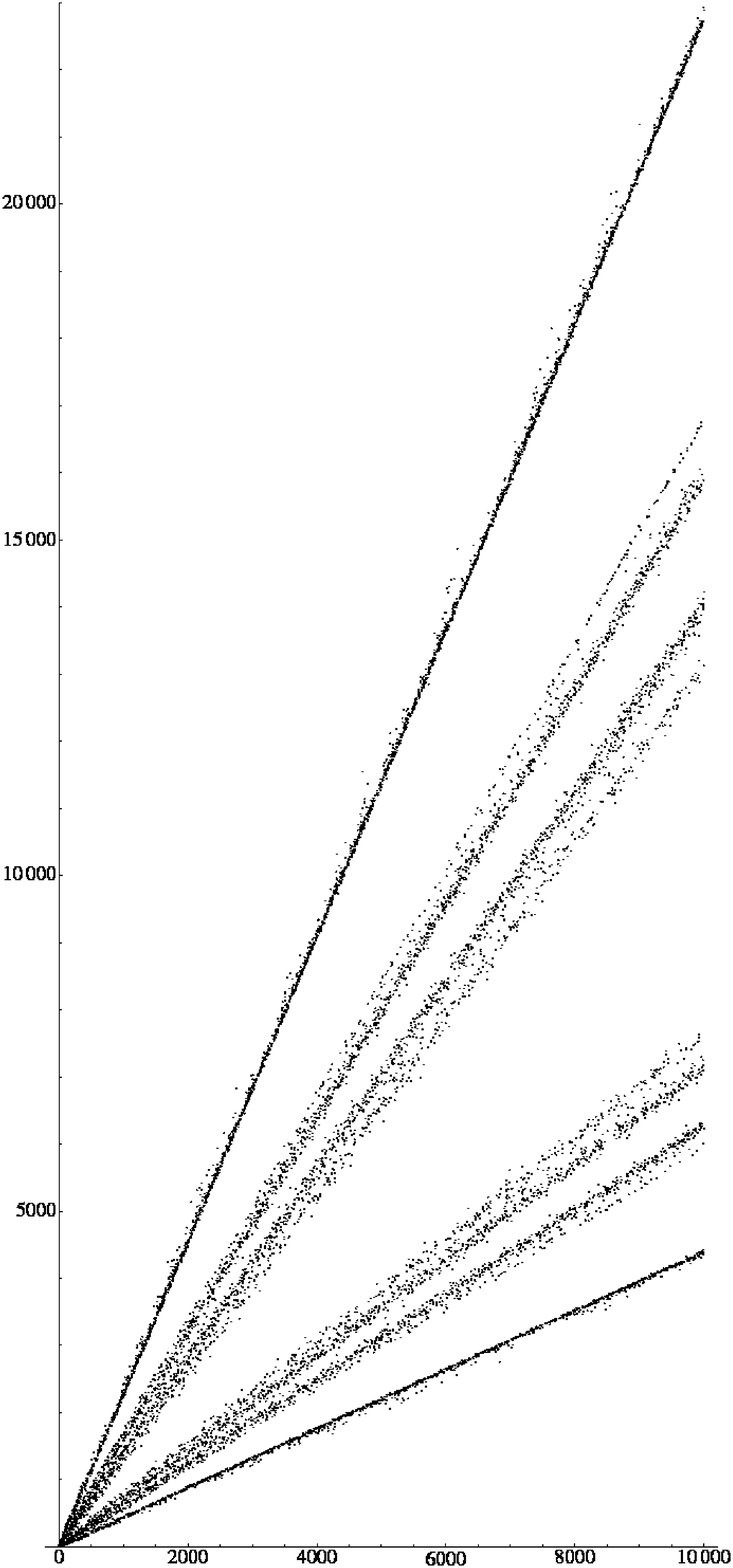}}  
\end{center}\caption{The $P$-positions  $(n,\pi_Q(n))$ 
of $\{(p,q)\mid p<q \le 5\}\G$ and $0\le n \le 10000$.}
\end{figure}
\clearpage
\begin{figure}[ht!]
\begin{center}
	  {\includegraphics[width=0.6\textwidth]{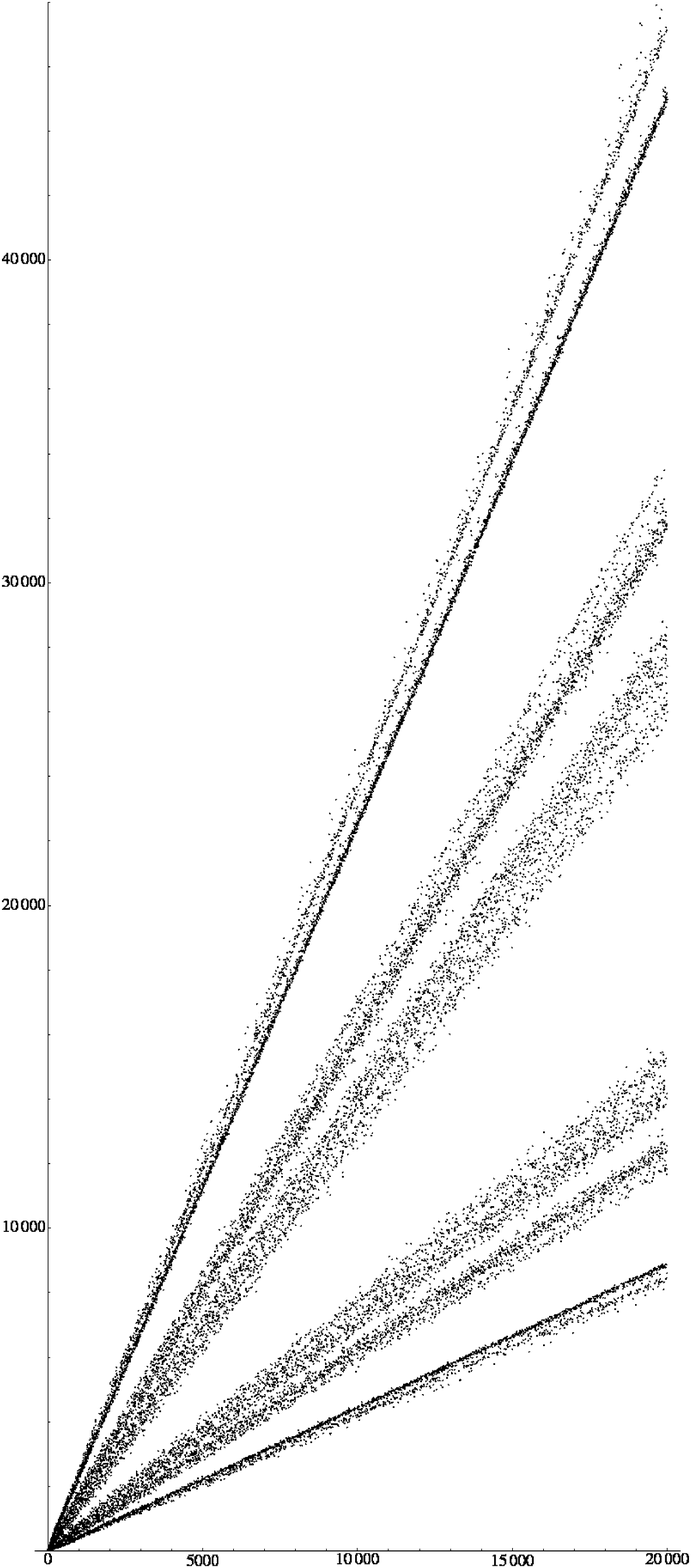}}  
\end{center}\caption{The $P$-positions  $(n,\pi_Q(n))$ 
of $\{(p,q)\mid p<q \le C\}\G$, $C=7$ and $0\le n \le 20000$. What happens 
for large values of $C$, will we get further splitting or will 
gradually the whole board become 'filled' with uniformly 
distributed $P$-positions?}
\end{figure}
\clearpage
\begin{figure}[ht!]
\begin{center}
	  {\includegraphics[width=0.7\textwidth]{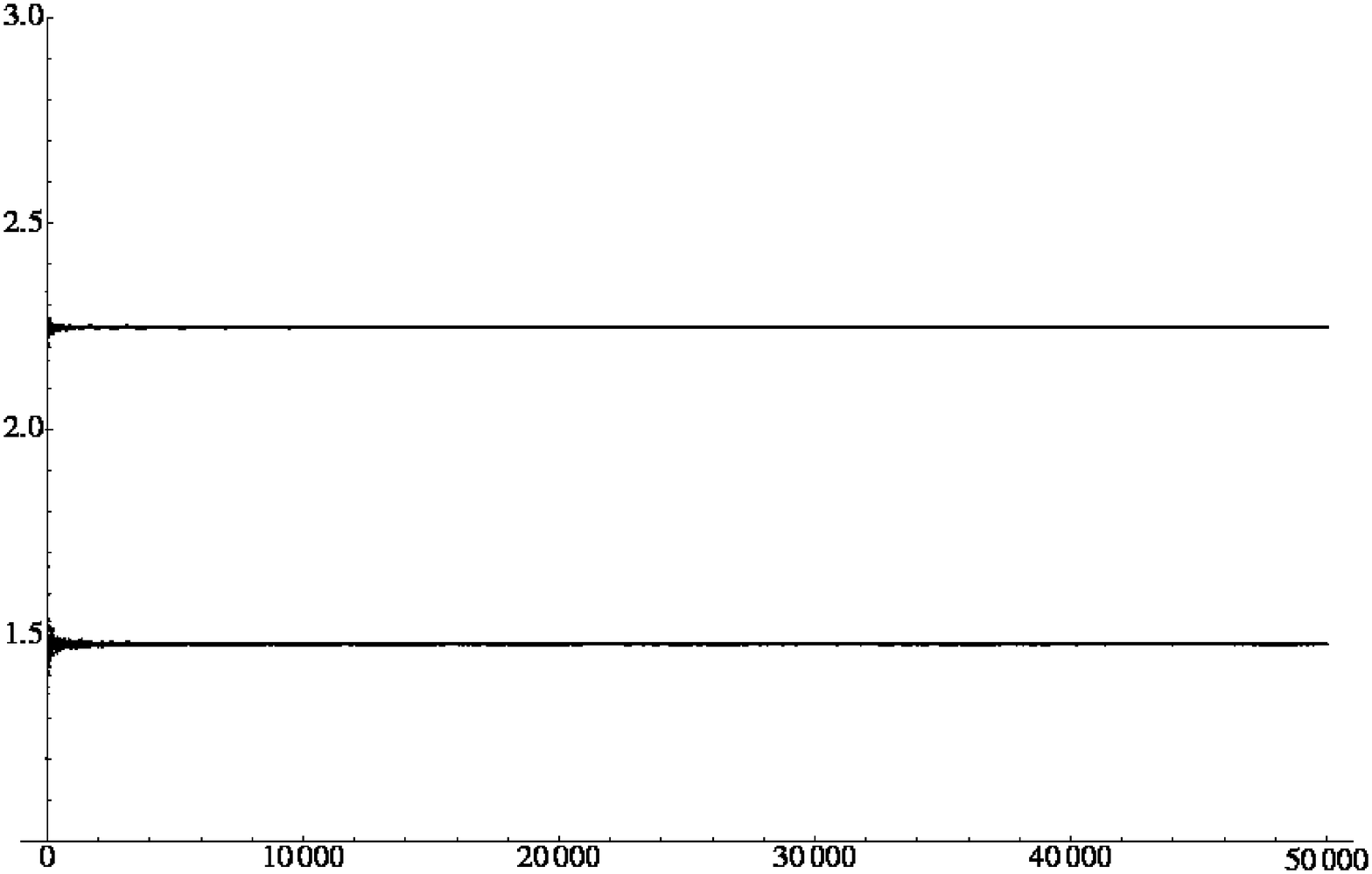}}  
\end{center}
\begin{center}
	  {\includegraphics[width=0.7\textwidth]{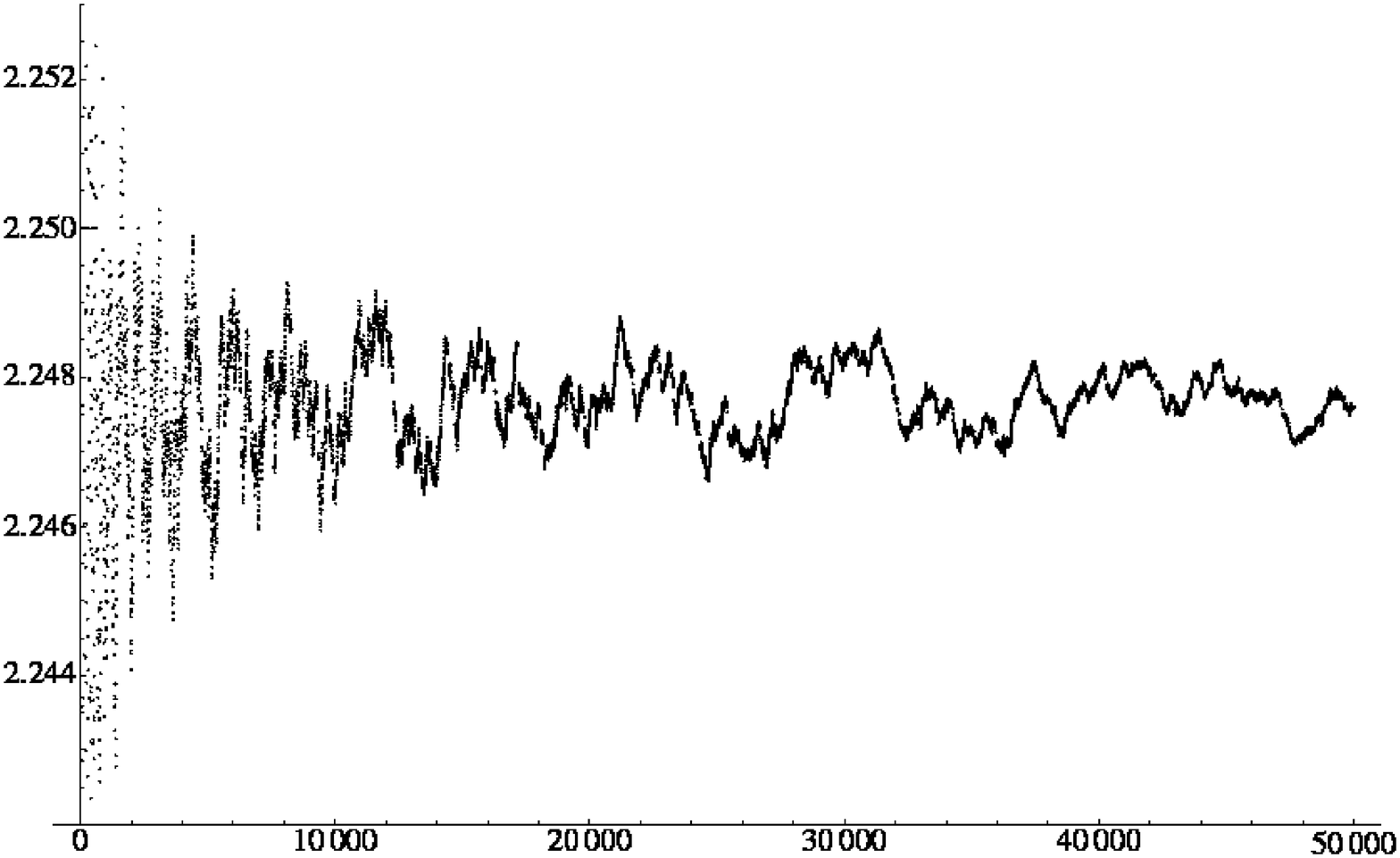}}  
\end{center}
\begin{center}
	  {\includegraphics[width=0.7\textwidth]{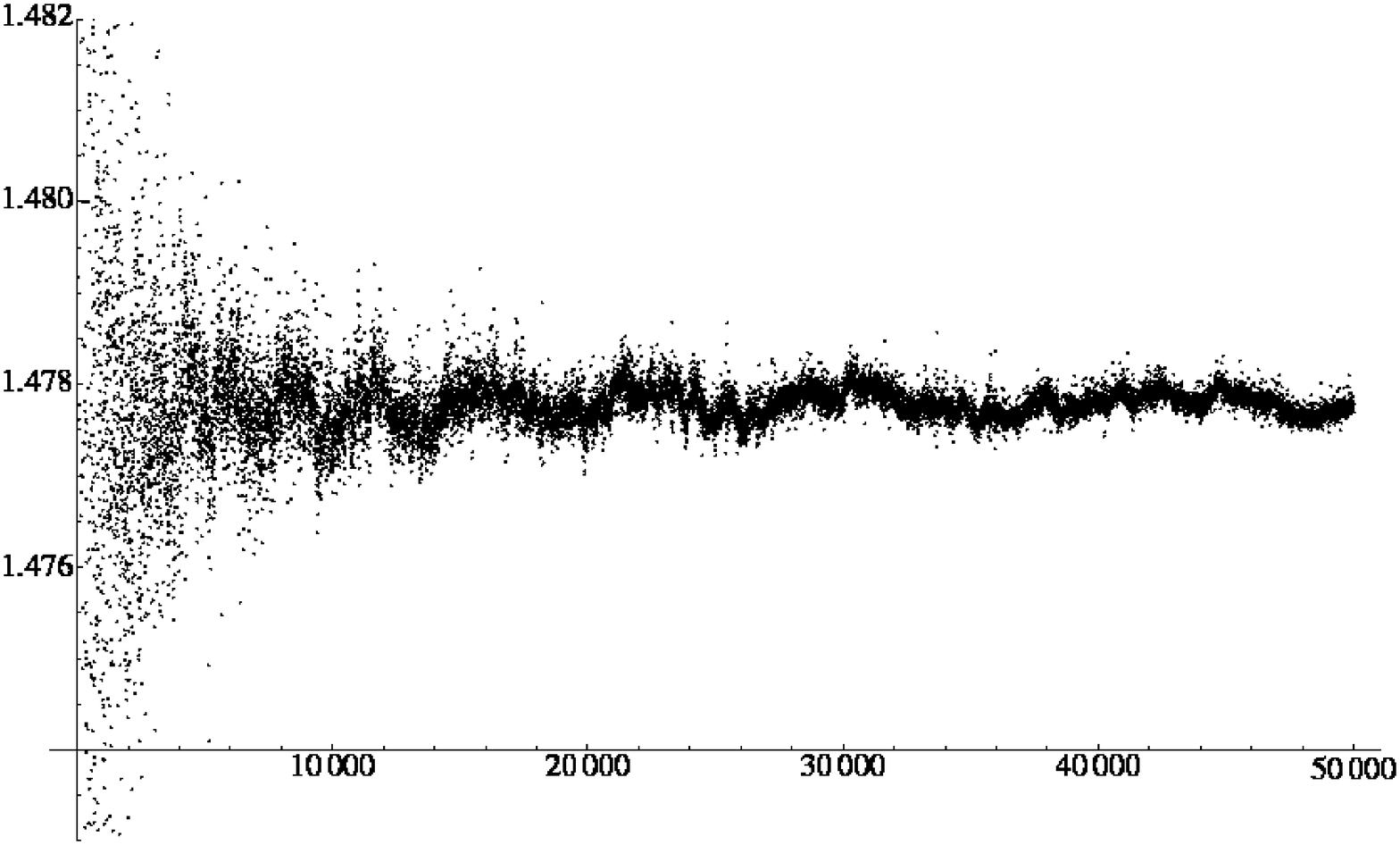}}  
\end{center}\caption{The figure at the top illustrates the 
ratio $b_n/a_n$ for $(1,2)\G$ and all $0\le n \le 50000$. We conjecture 
that there exist two complementary sequences $u$ (middle figure) and $l$ 
(lower) such that  $b_{u_i}/a_{u_i}\rightarrow 2.247\ldots $ 
(roughly 40\%) 
and $\lim_{i\rightarrow \infty }b_{l_i}/a_{l_i}\rightarrow 1.478\ldots$ 
(roughly 60\%).}\label{A12}
\end{figure}
\clearpage
\begin{figure}[ht!]
\begin{center}
	  {\includegraphics[width=0.6\textwidth]{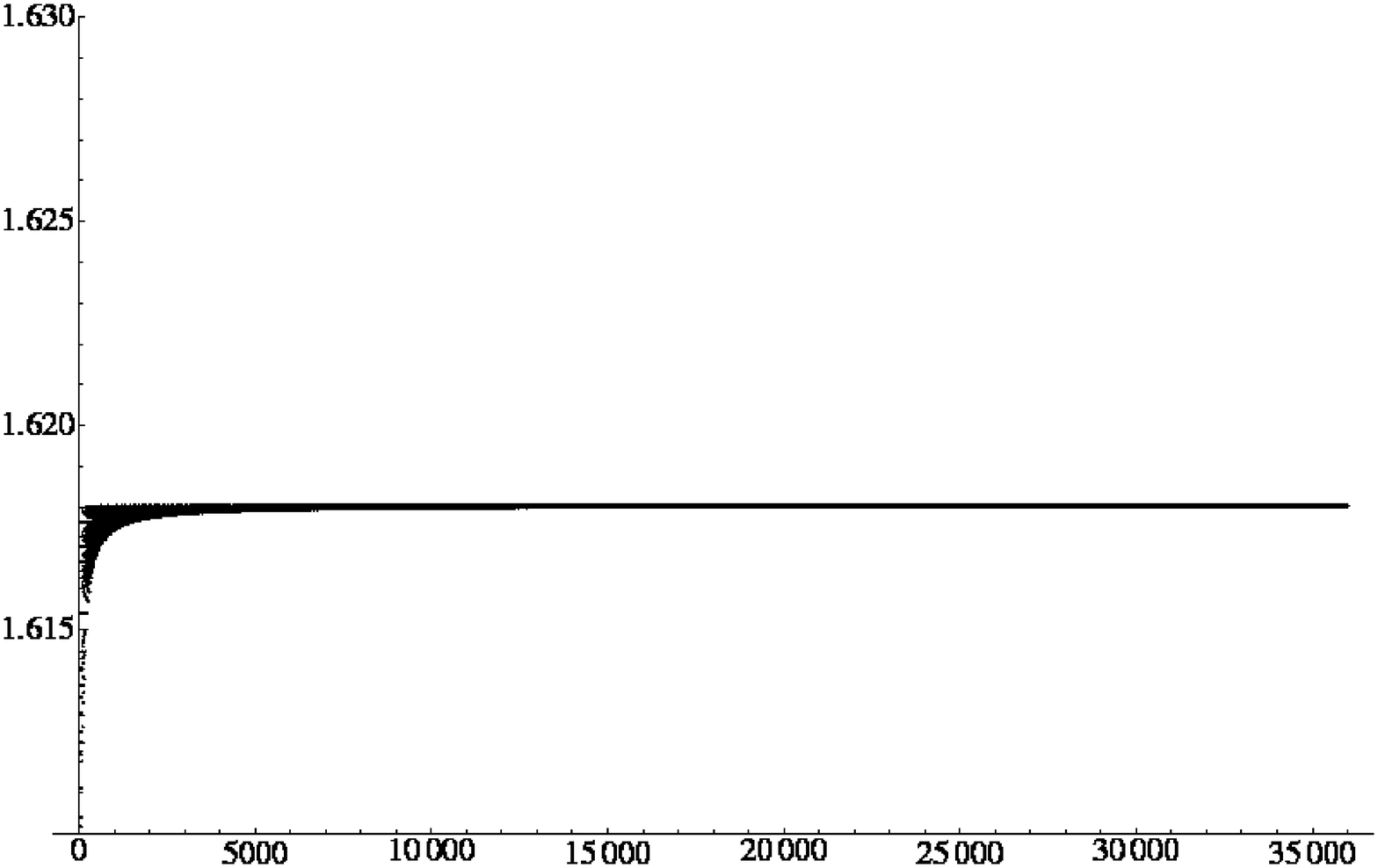}}  
\end{center}
\begin{center}
	  {\includegraphics[width=0.6\textwidth]{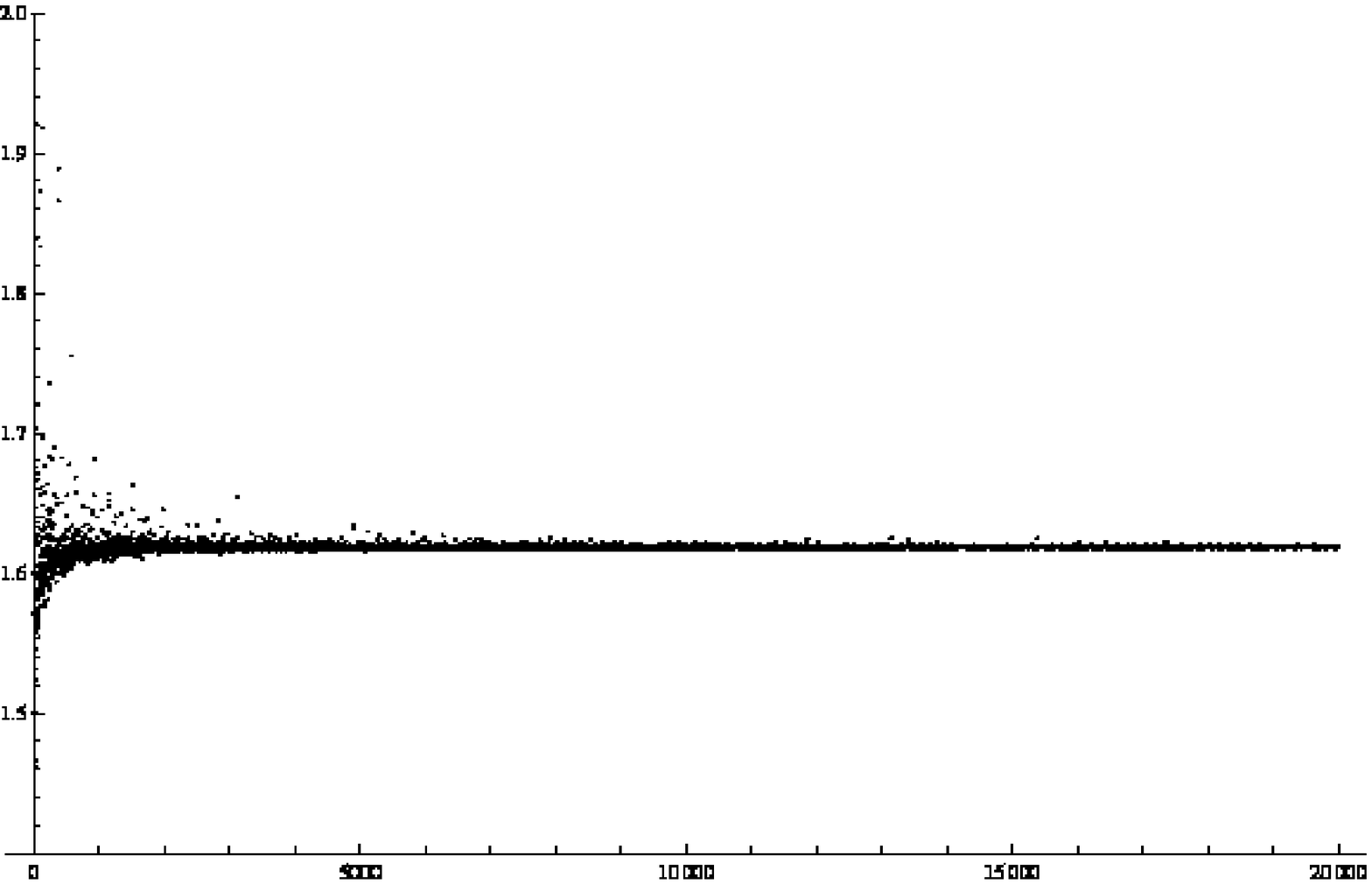}}  
\end{center}
\begin{center}
	  {\includegraphics[width=0.6\textwidth]{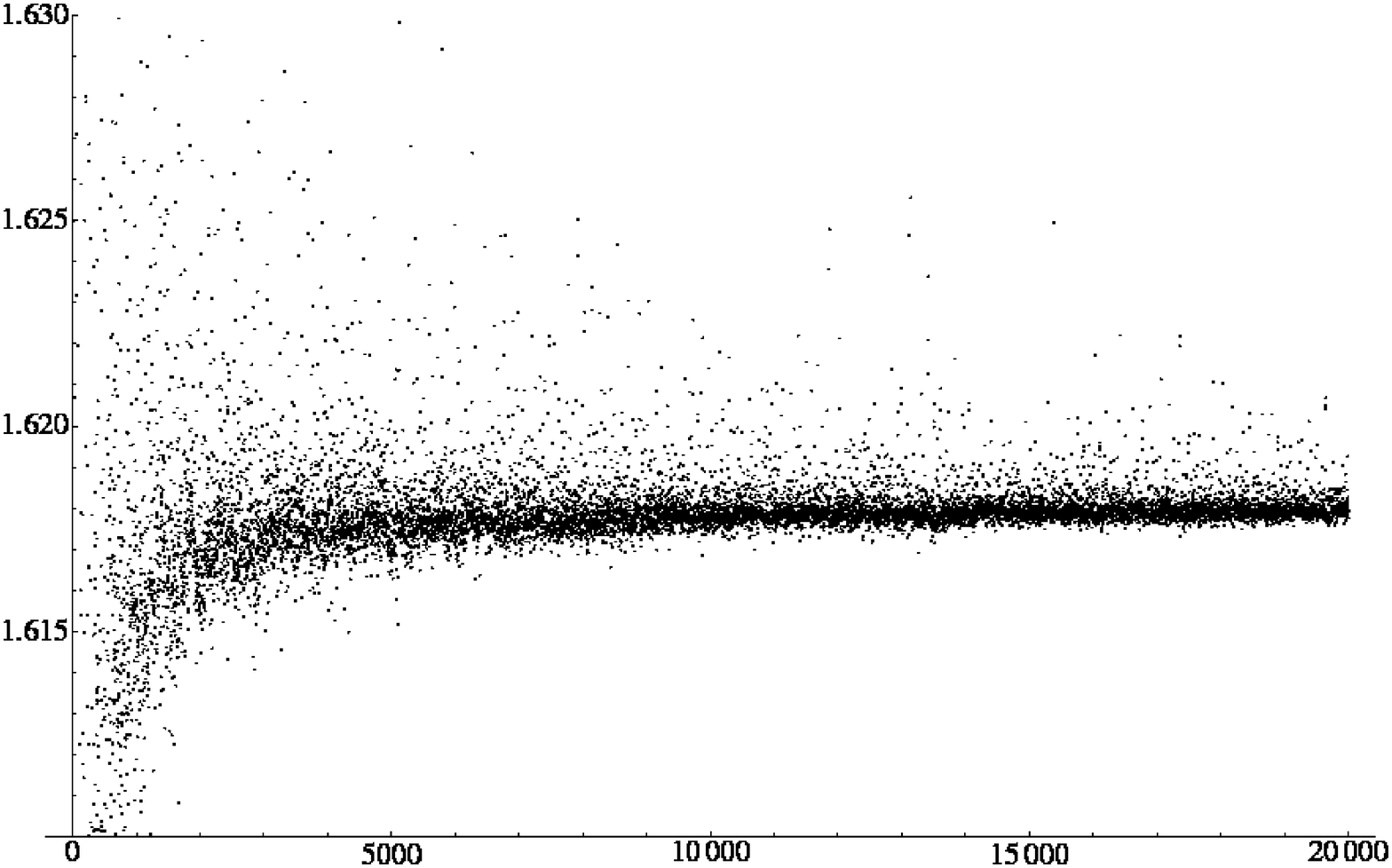}}  
\end{center}\caption{The figure at the top illustrates the ratio $b_n/a_n$ 
for Wythoff Nim and all equivalent games $(p,q)\G$, that is whenever $(p,q)$ 
is a non-splitting pair and $q/p < \Phi$. The two lower figures illustrate 
the corresponding ratios for $(2,4)\G$ 
and all $0\le n \le 20000$. Our simulation suggests no 
split, rather $b_{n}/a_{n}\rightarrow 1.618\ldots$. However, the latter 
game is clearly not equivalent to Wythoff Nim. 
See also Table \ref{A3}.}\label{A13}
\end{figure}
\clearpage

\begin{figure}[ht!]
\begin{center}
	  {\includegraphics[width=0.7\textwidth]{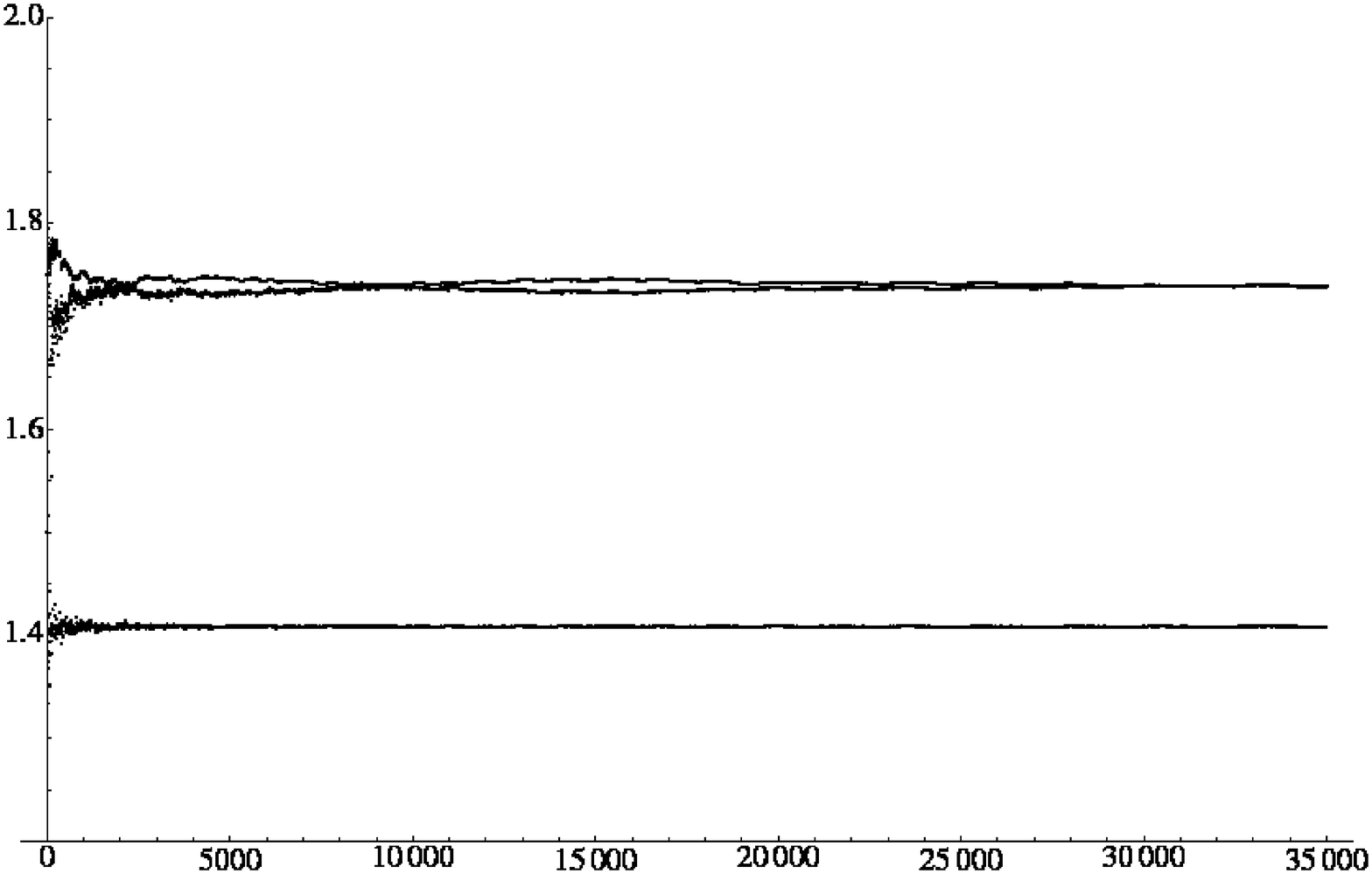}}  
\end{center}
\begin{center}
	  {\includegraphics[width=0.7\textwidth]{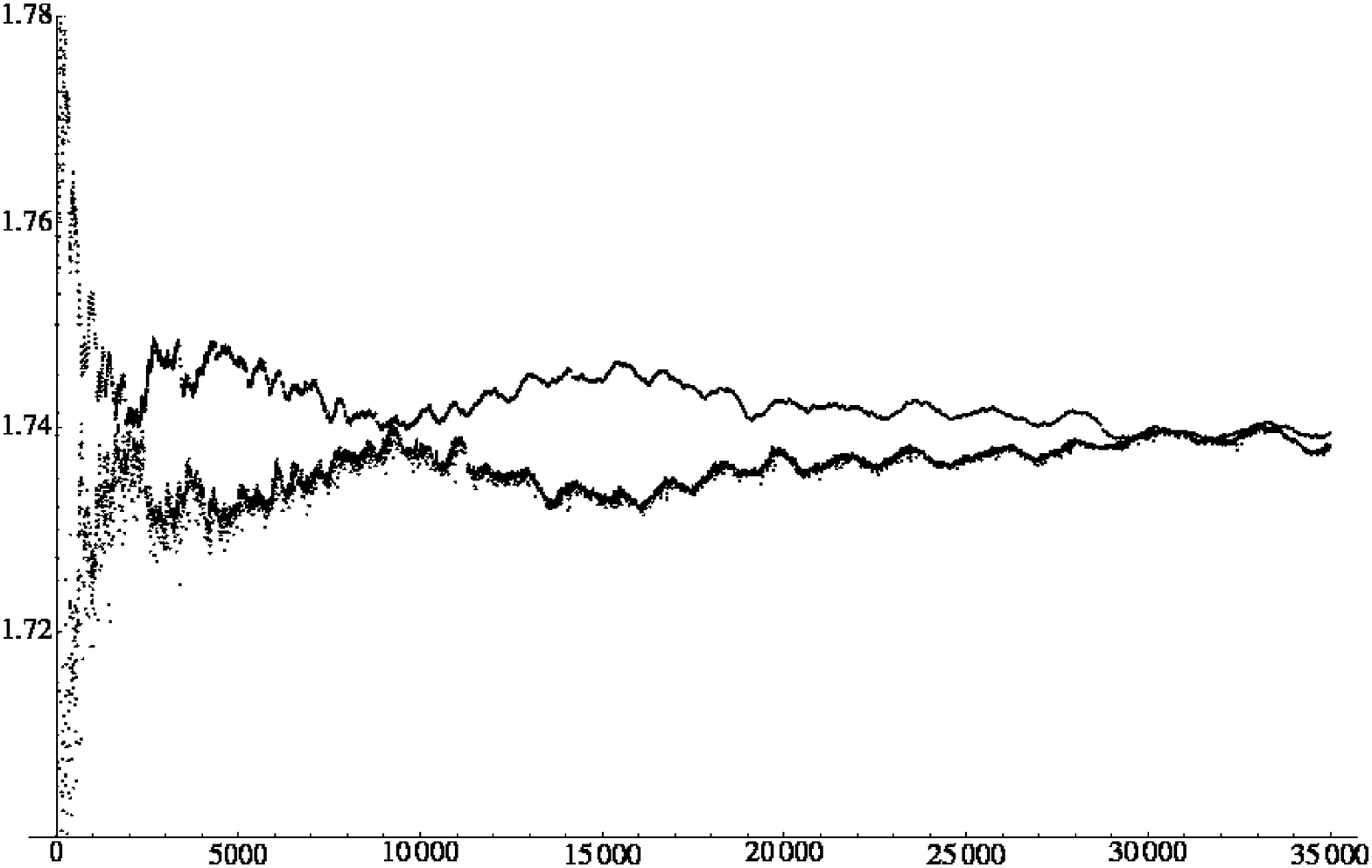}}  
\end{center}
\begin{center}
	  {\includegraphics[width=0.7\textwidth]{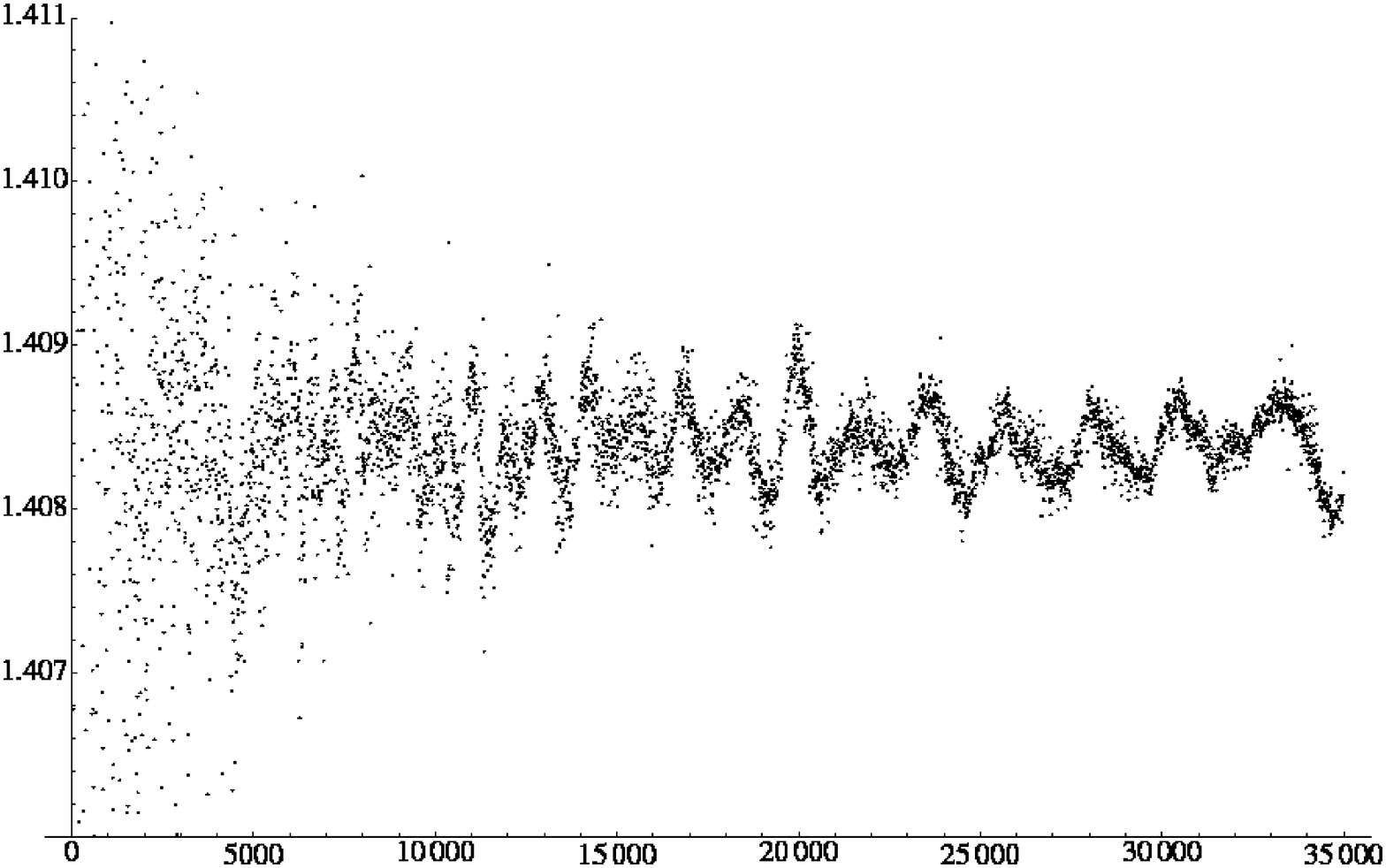}}  
\end{center}\caption{The ratio $b_n/a_n$ for $(2,3)\G$ and 
 $0\le a_n \le 35000$. We conjecture that there is a pair of 
complementary sequences $u$ (middle picture) and $l$ (lower picture) 
such that $b_{u_i}/a_{u_i}\rightarrow 1.74\ldots $ (roughly 80\%) and
$b_{l_i}/a_{l_i}\rightarrow 1.408\ldots $ (roughly 20\%).}
\end{figure}
\clearpage
\begin{figure}
\begin{center}
	  {\includegraphics[width=0.8\textwidth]{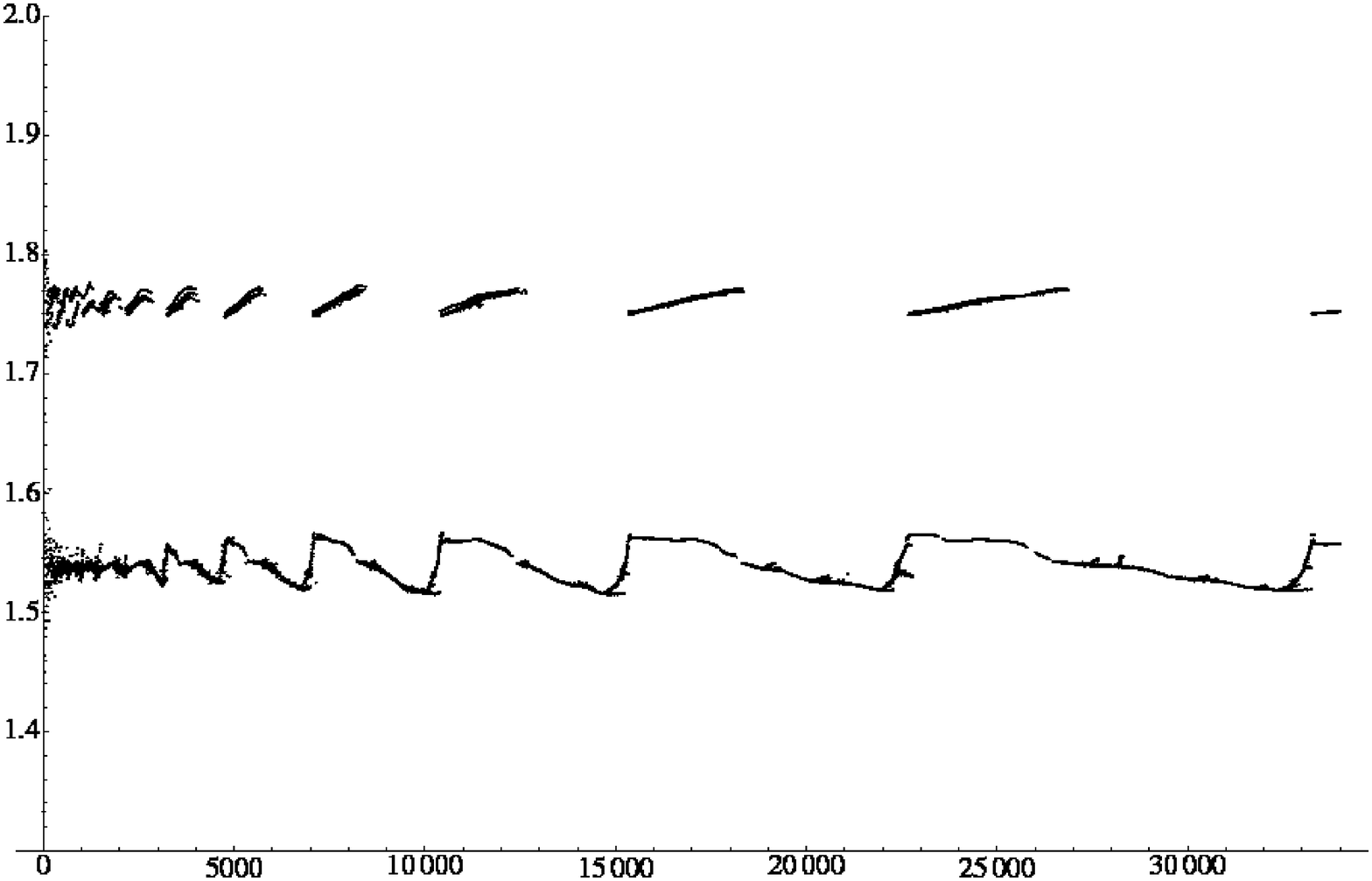}}  
\end{center}\caption{The ratio $b_n/a_n$ for $(3,5)\G$ and 
$0\le a_n \le 35000$.}
\end{figure}

\begin{figure}[ht!]
\begin{center}
	  {\includegraphics[width=0.8\textwidth]{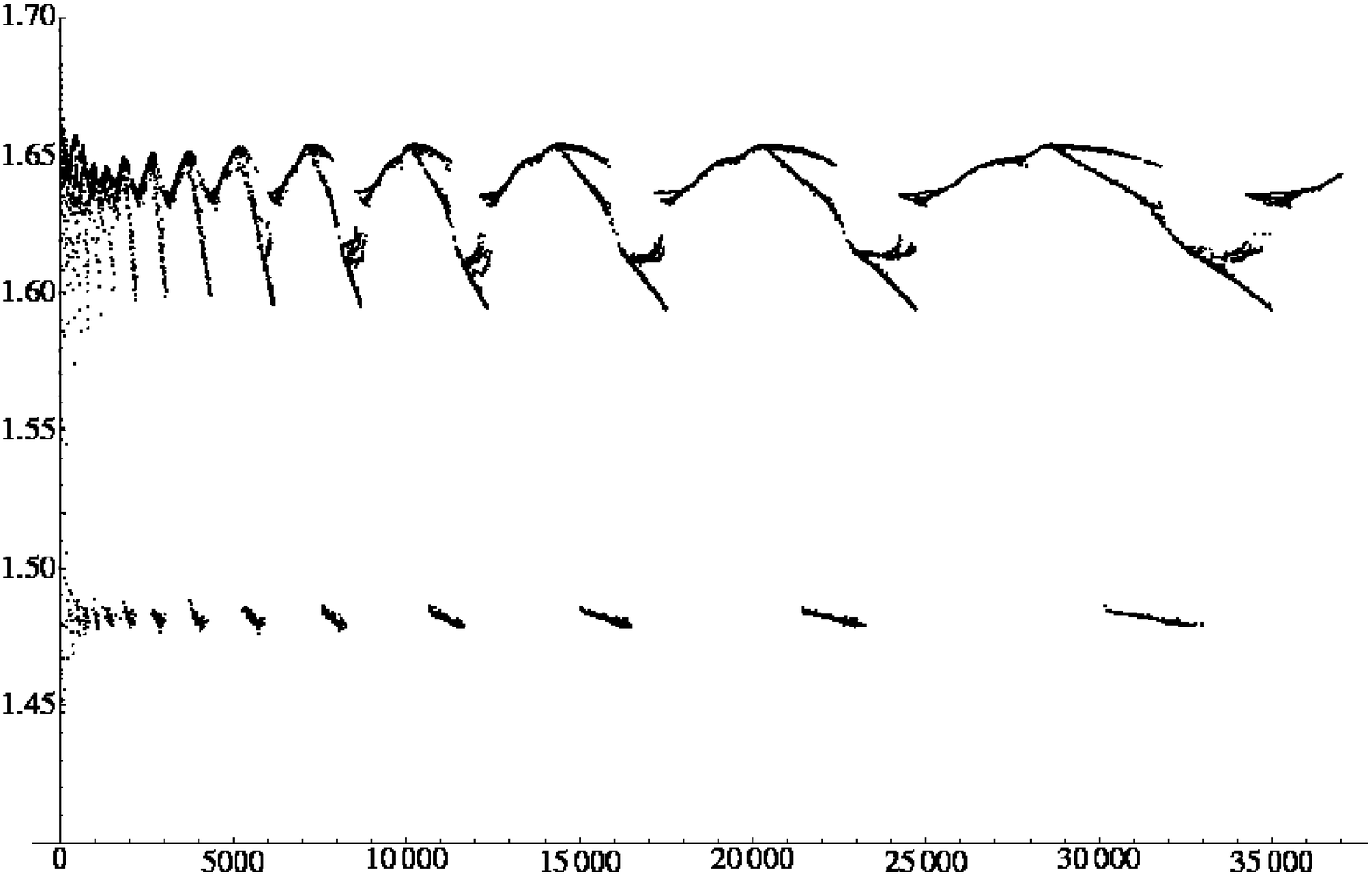}}  
\end{center}\caption{The ratio $b_n/a_n$ for $(4,6)\G$ 
and $0\le a_n \le 35000$. Our data seems to suggest  that there is a pair 
of complementary sequences 
$u$ and $l$ such that for large $i$, 
$$ 1.60\ldots < b_{u_i}/a_{u_i} < 1.66\ldots $$ and the quotient 
is 'drifting back and forth' in this interval, but 
$b_{l_i}/a_{l_i}\rightarrow 1.48\ldots $ as $i\rightarrow \infty$.}
\end{figure}
\clearpage
\begin{figure}
\begin{center}
	  {\includegraphics[width=0.8\textwidth]{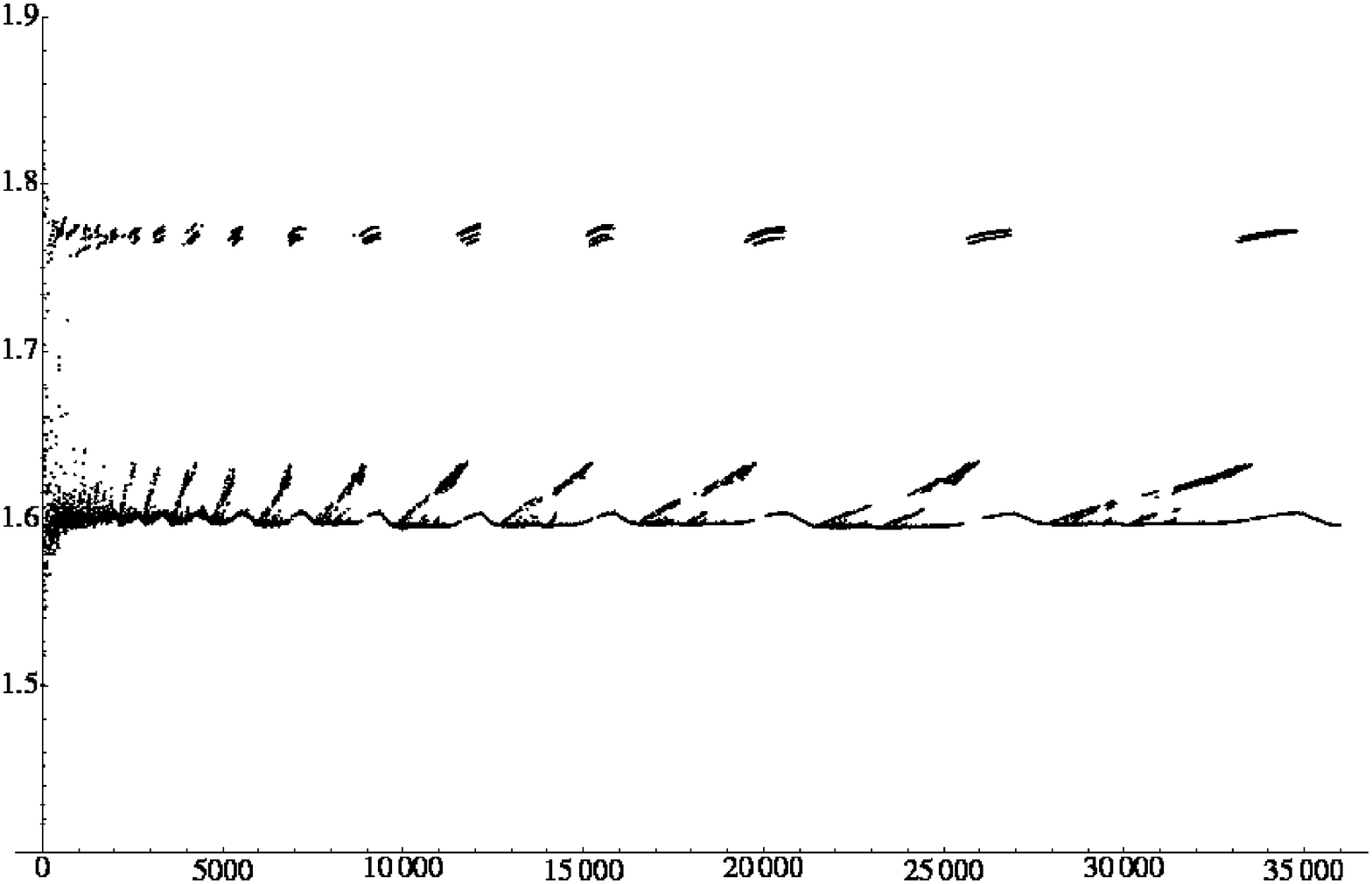}}  
\end{center}\caption{The ratio $b_n/a_n$ for $(4,7)\G$ and $0\le n \le 35000$. 
It does seem to split asymptotically, 
but maybe only the weaker form of our conjecture holds for this case, namely 
our data suggest that there is a pair of complementary sequences 
$u$ and $l$ such that $b_{l_i}/a_{l_i}$ is 'drifting' in 
the interval $[1.59, 1.63]$ for 'large' $i$, 
but $ b_{u_i}/a_{u_i}\rightarrow 1.77\ldots $ .}
\end{figure}

\begin{figure}[ht!]
\begin{center}
	  {\includegraphics[width=0.8\textwidth]{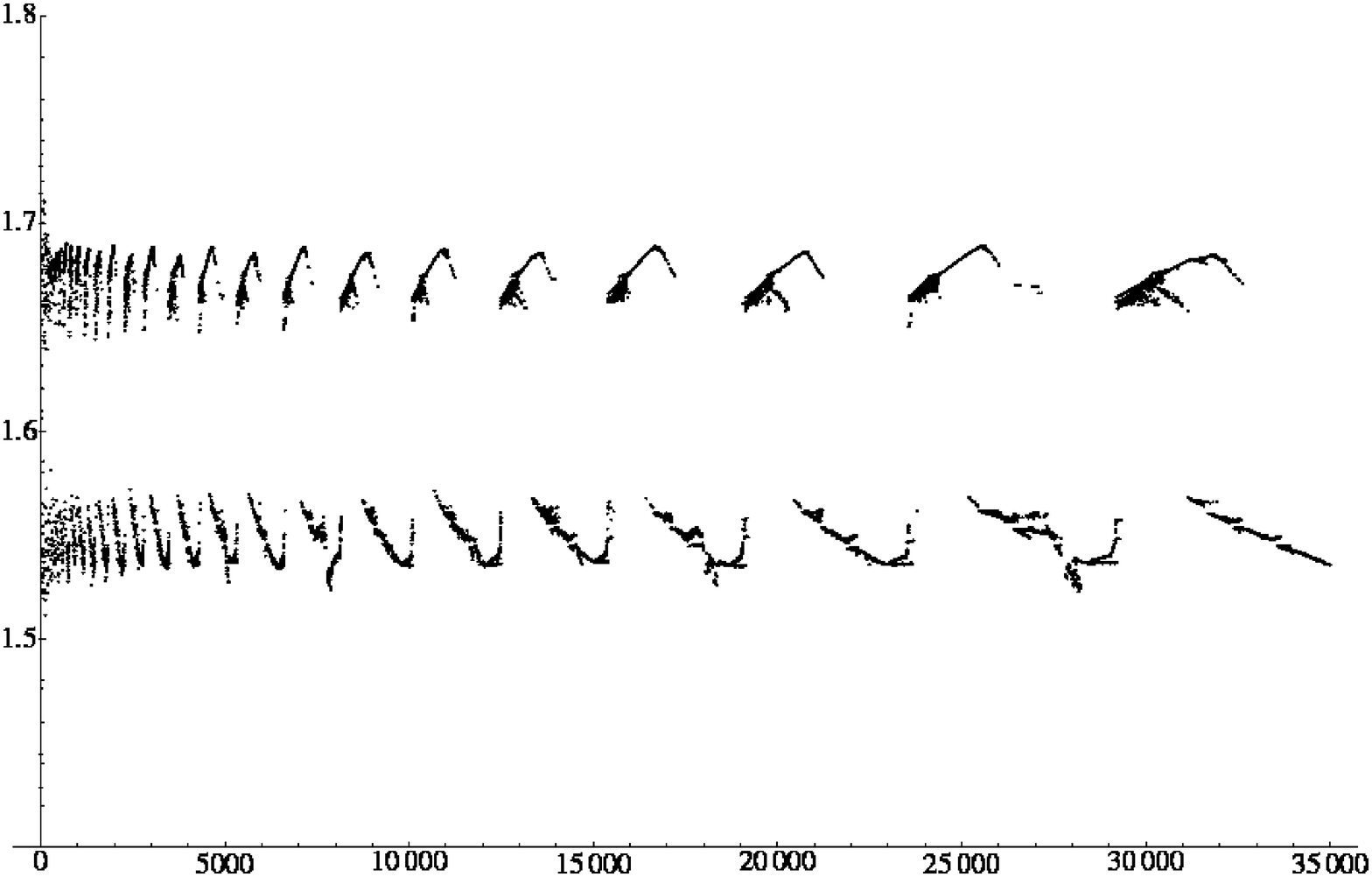}}  
\end{center}\caption{The ratio $b_n/a_n$ for $(5,8)\G$ and $0\le n \le 35000$.}
\end{figure}
\clearpage
\begin{figure}
\begin{center}
	  {\includegraphics[width=0.8\textwidth]{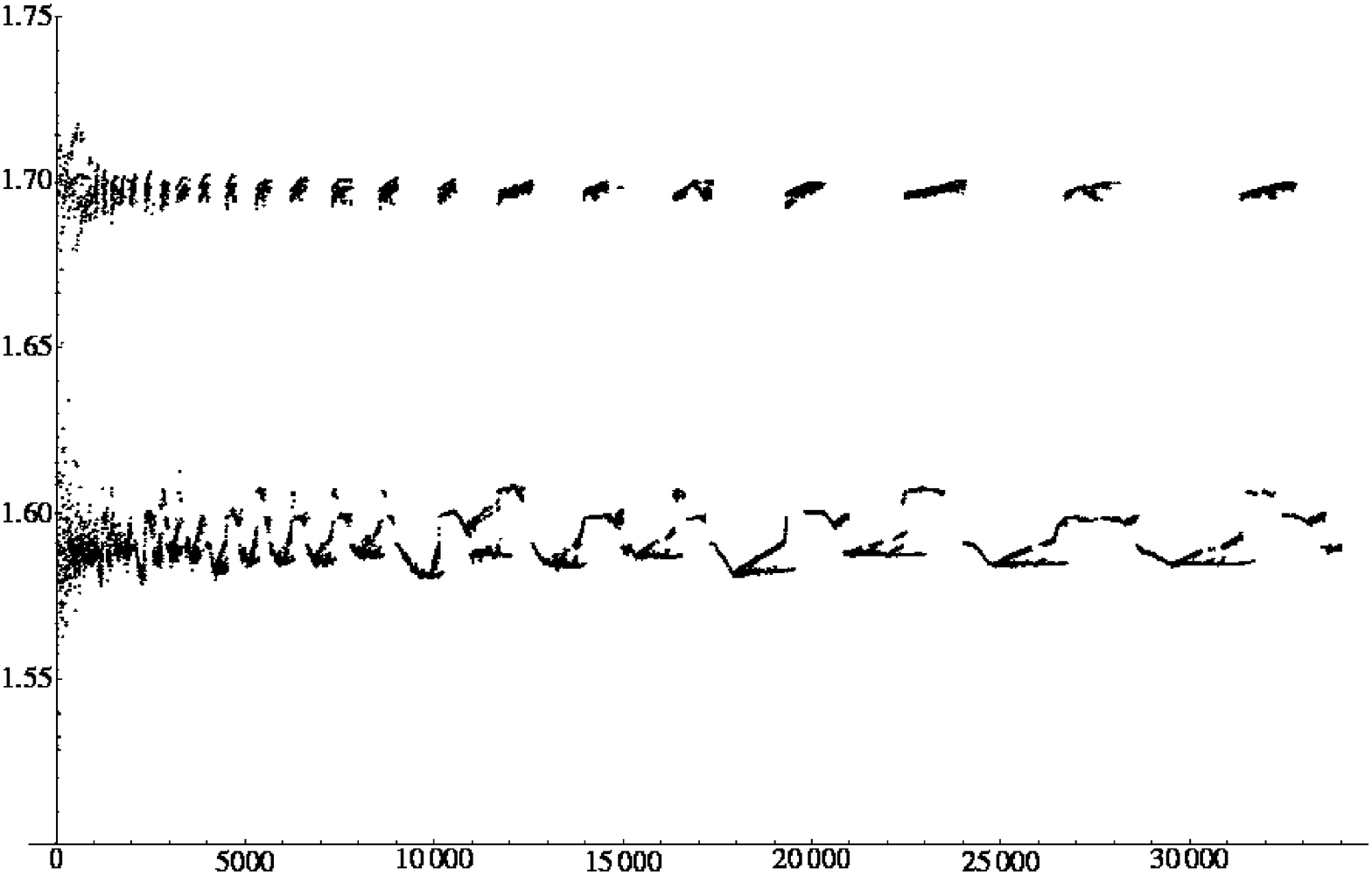}}  
\end{center}\caption{The ratio $b_n/a_n$ for $(6,10)\G$ and $0\le n \le 35000$.}
\end{figure}

\begin{figure}[ht!]
\begin{center}
	  {\includegraphics[width=0.8\textwidth]{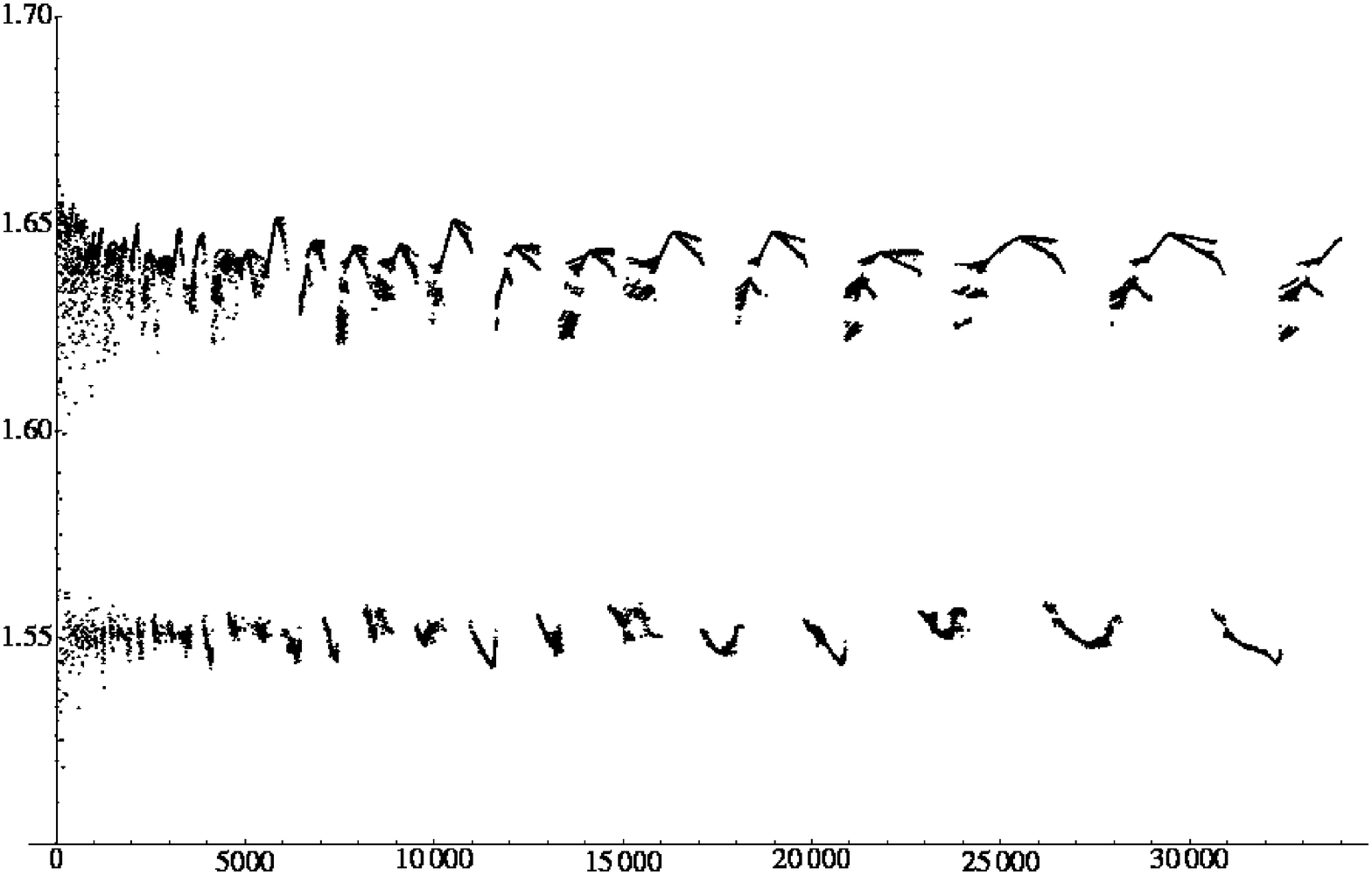}}  
\end{center}\caption{The ratio $b_n/a_n$ for $(7,11)\G$ and $0\le n \le 35000$.}
\end{figure}
\clearpage
\begin{figure}
\begin{center}
	  {\includegraphics[width=0.8\textwidth]{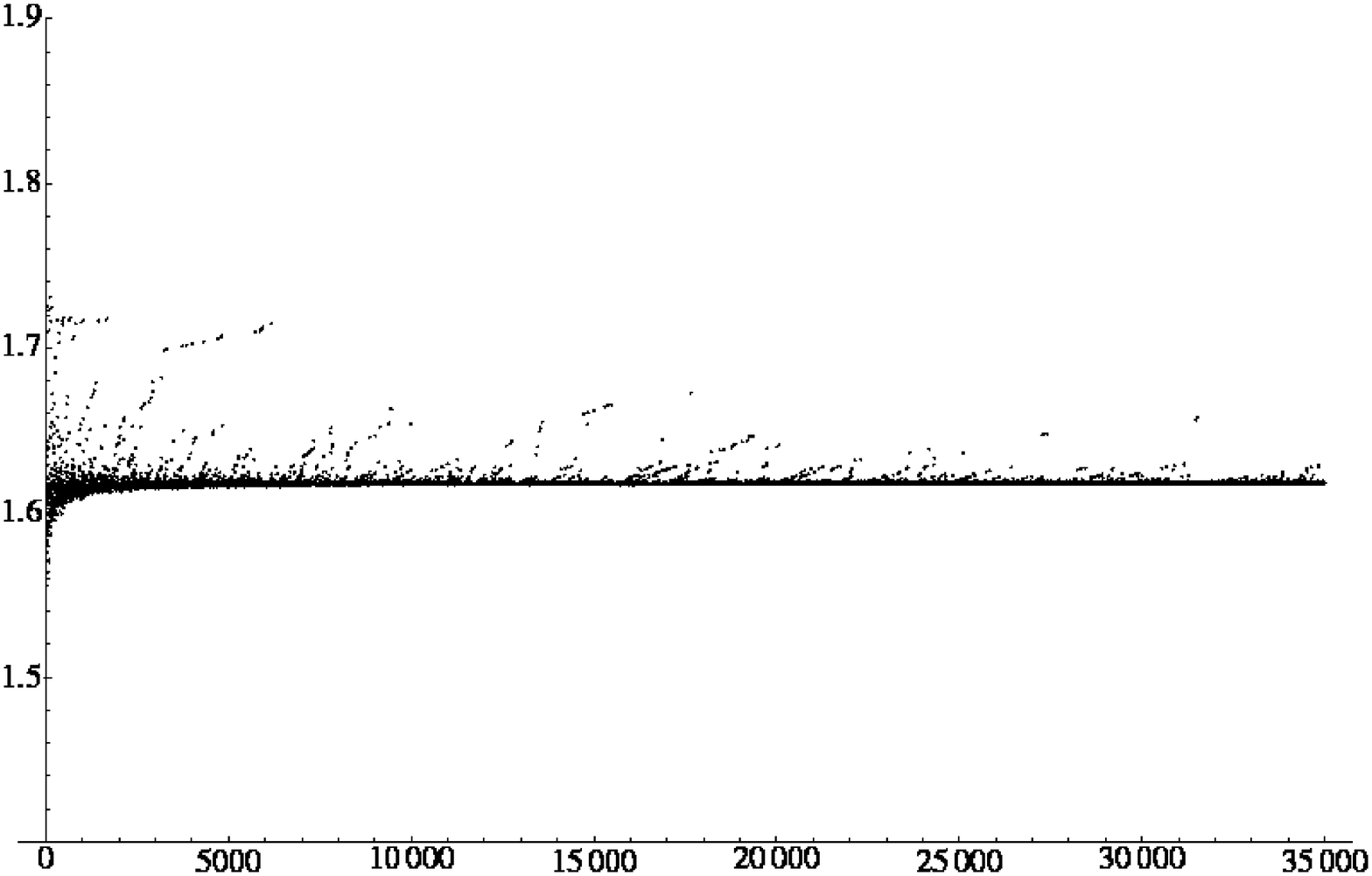}}  
\end{center}
\begin{center}
	  {\includegraphics[width=0.8\textwidth]{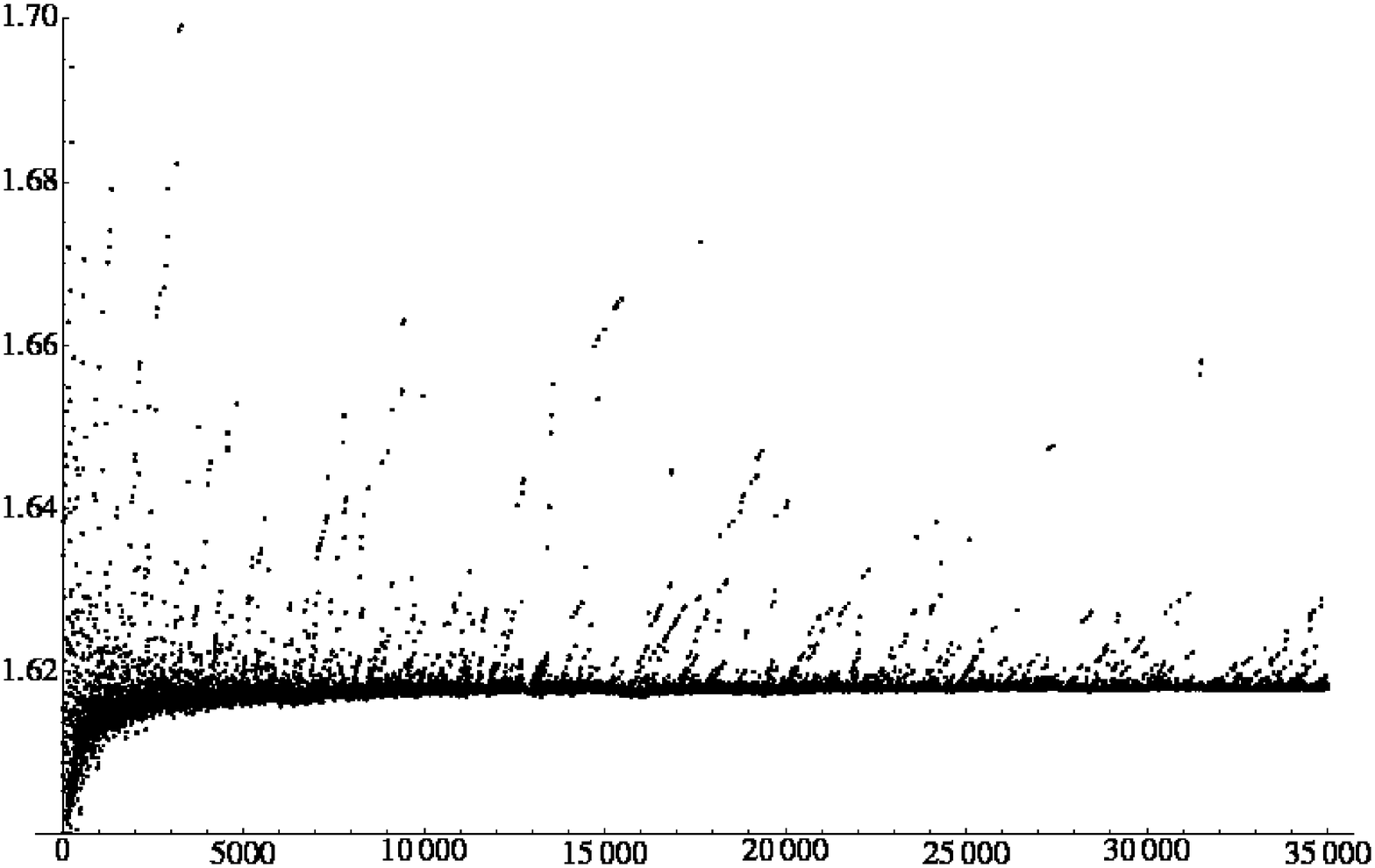}}  
\end{center}\caption{The ratio $b_n/a_n$ for $(7,12)\G$ 
and $0\le n \le 35000$.}\label{A21}
\end{figure}
\clearpage
\begin{figure}[ht!]
\begin{center}
	  {\includegraphics[width=0.8\textwidth]{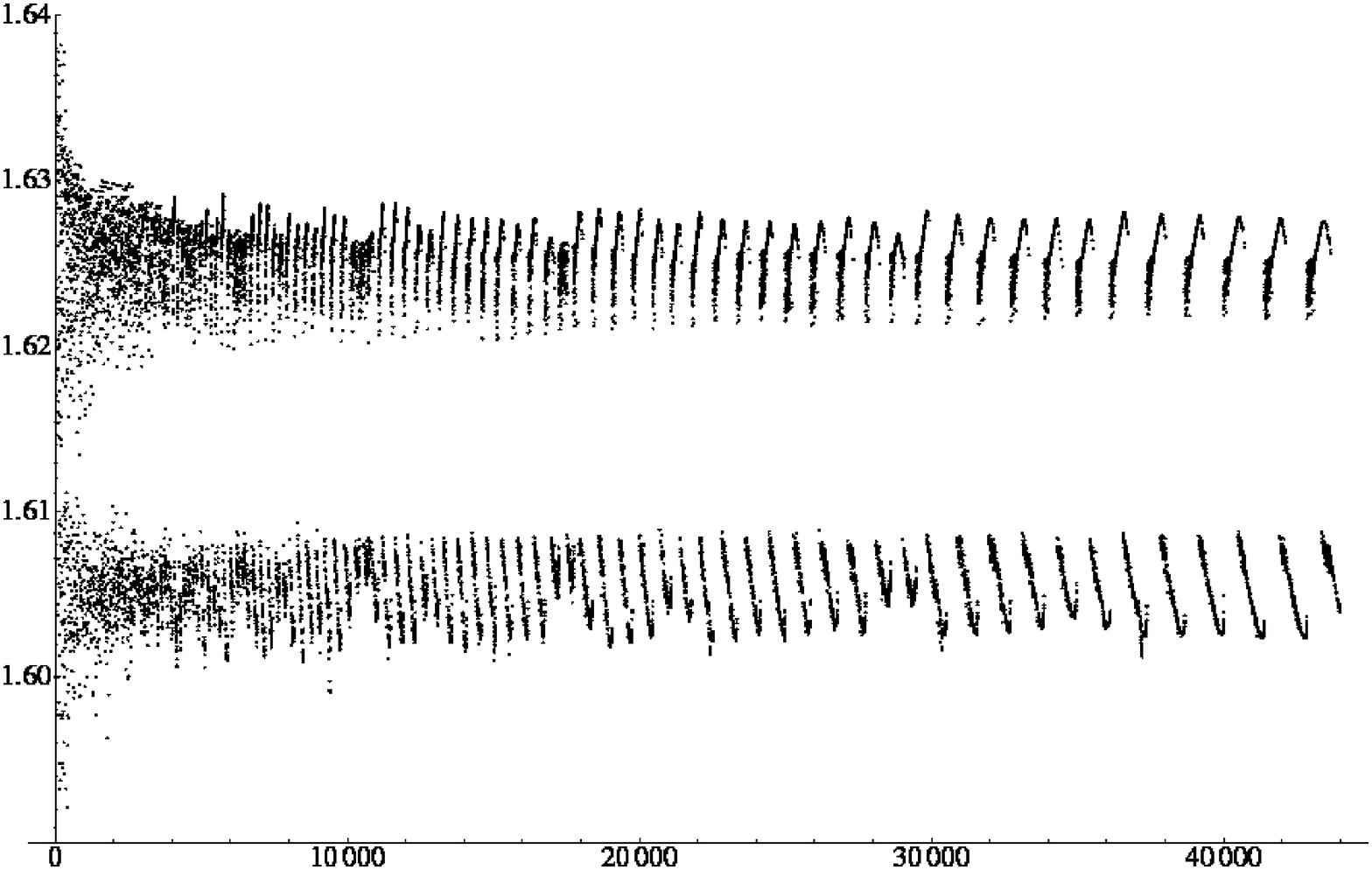}}  
\end{center}\caption{The ratio $b_n/a_n$ for $(31,50)\G$ and $0\le n \le 50000$.}
\end{figure}

\begin{figure}[ht!]
\begin{center}
	  {\includegraphics[width=0.8\textwidth]{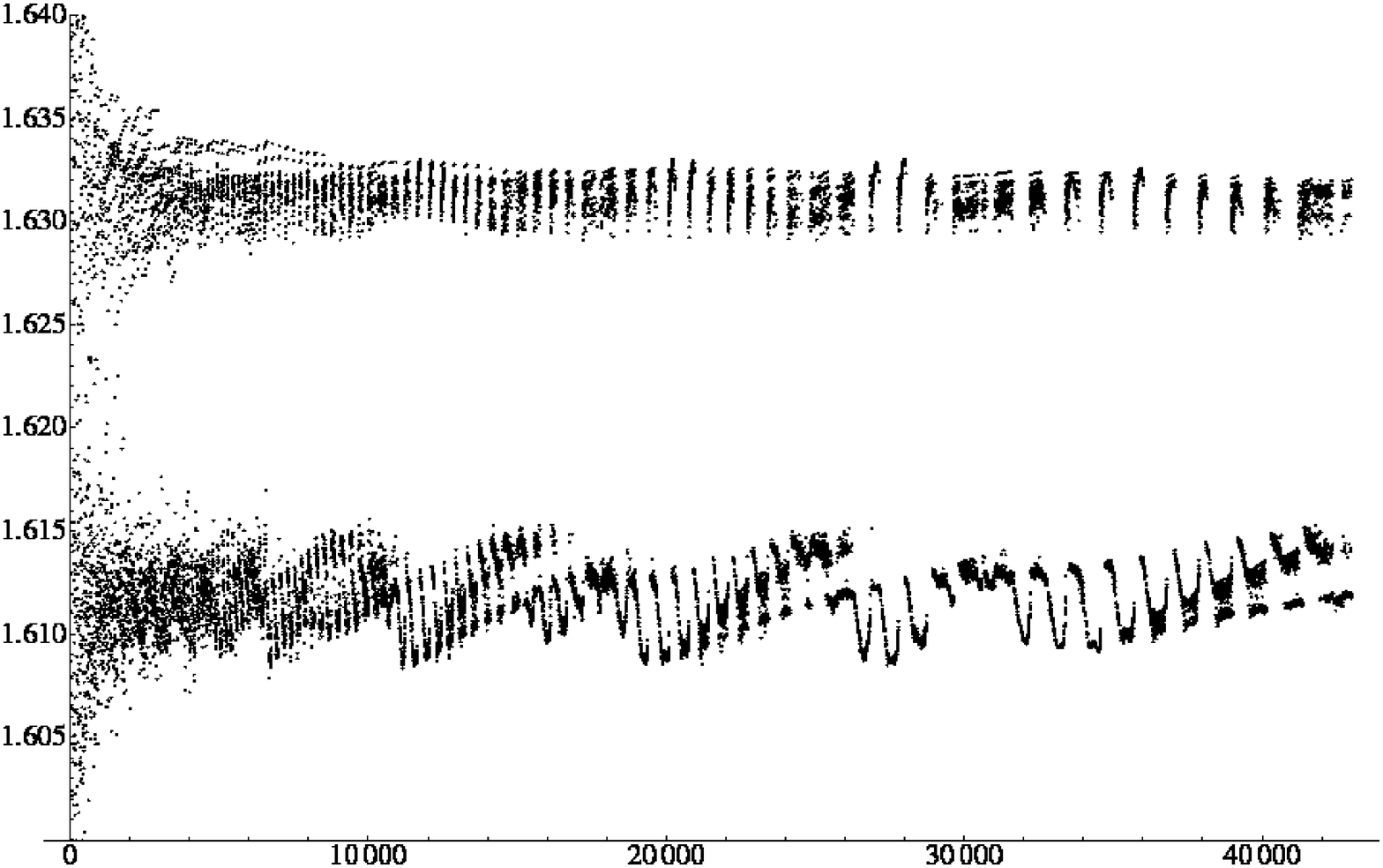}}  
\end{center}\caption{The ratio $b_n/a_n$ for $(32,52)\G$ and 
$0\le n \le 50000$.}
\end{figure}
\clearpage
\begin{figure}
\begin{center}
	  {\includegraphics[width=0.8\textwidth]{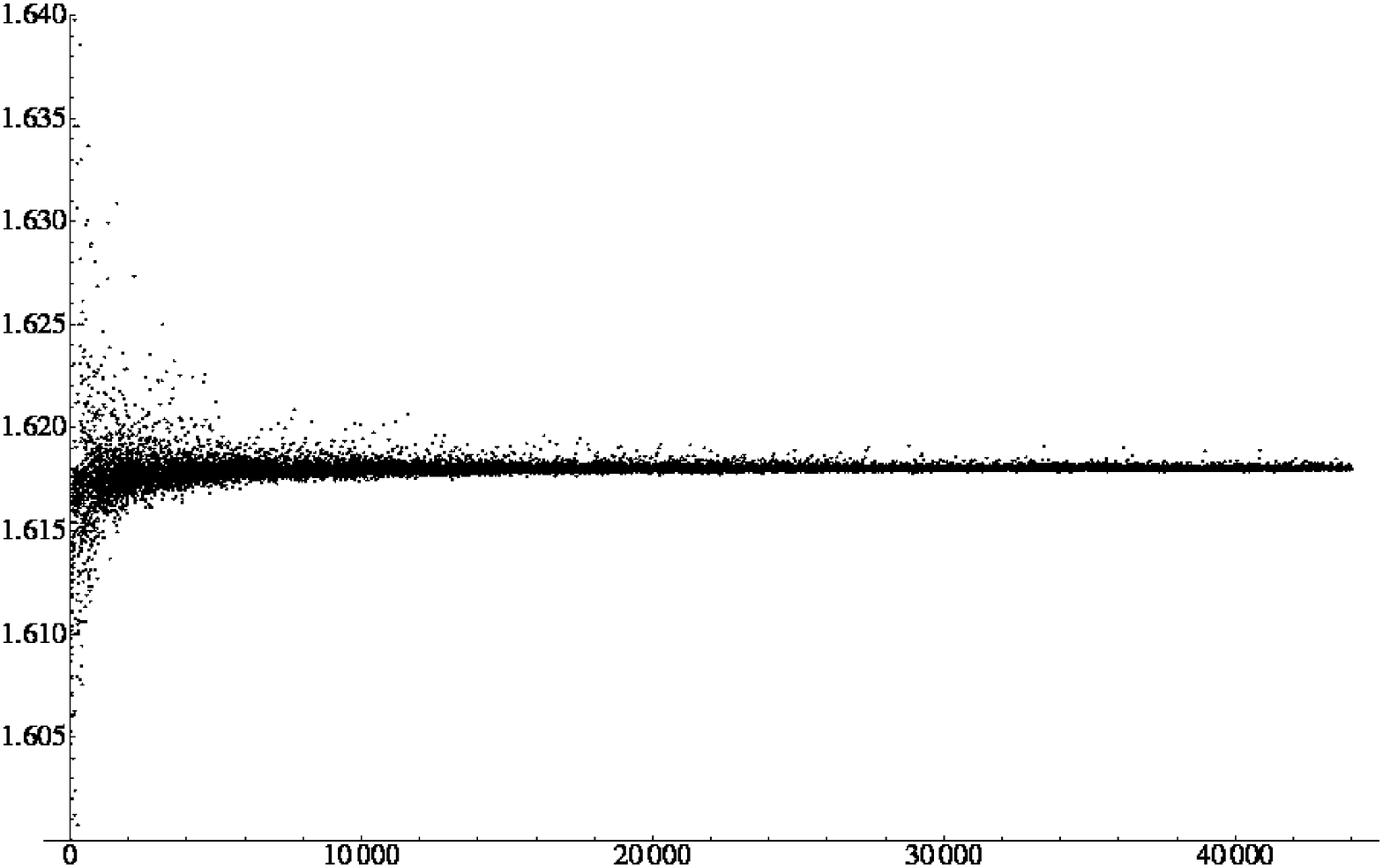}}  
\end{center}
\begin{center}
	  {\includegraphics[width=0.8\textwidth]{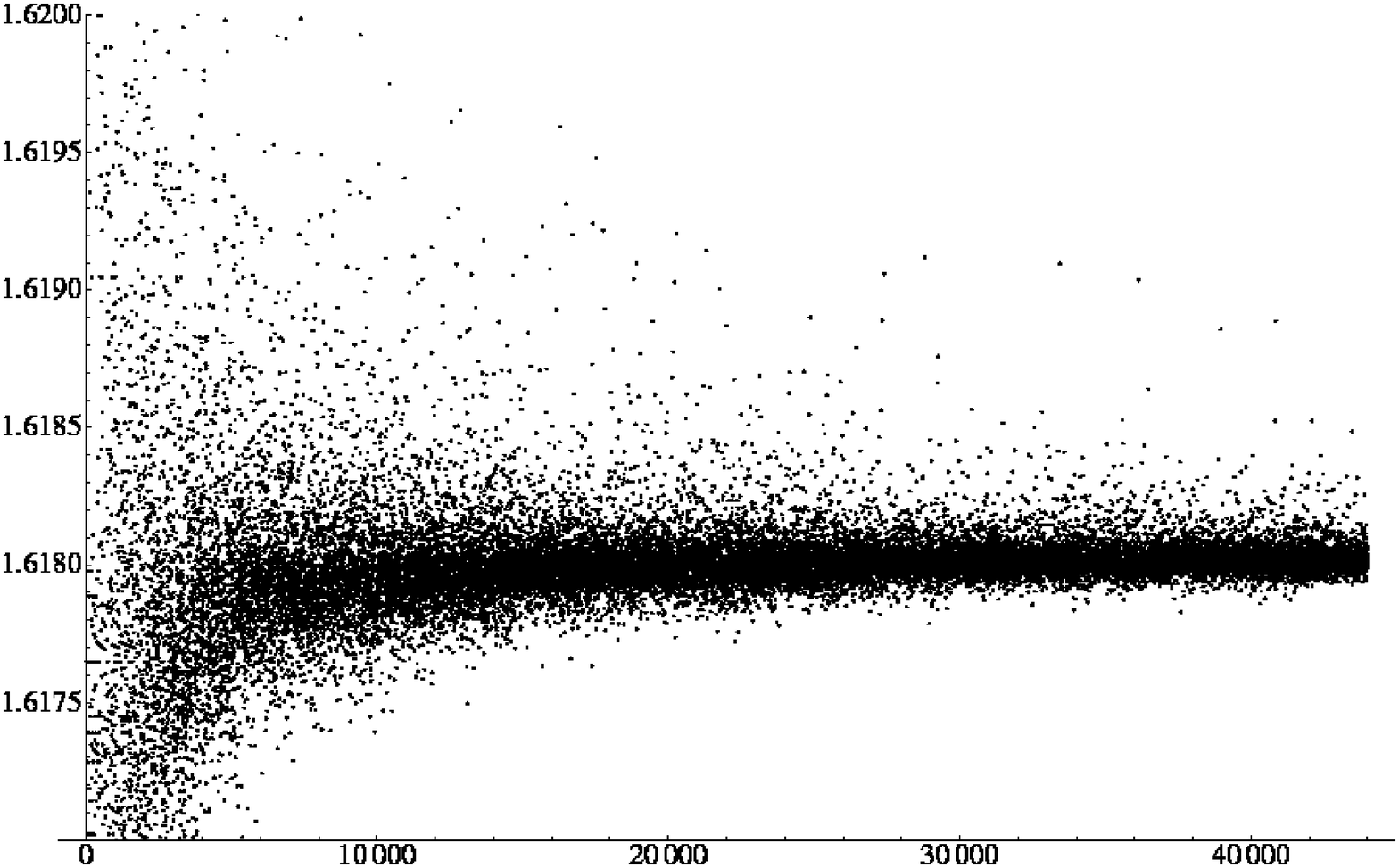}}  
\end{center}\caption{The ratio $b_n/a_n$ for $(31,51)\G$ 
and $0\le n \le 50000$. Notice that $(31,51)$ is a non-splitting pair, 
but $51/31 > 1.645 > \Phi$. As in Figures \ref{A13} and \ref{A21}, one 
may observe some perturbation of the $P$-positions of Wythoff Nim.}
\end{figure}
\clearpage
\begin{figure}
\begin{center}
	  {\includegraphics[width=0.8\textwidth]{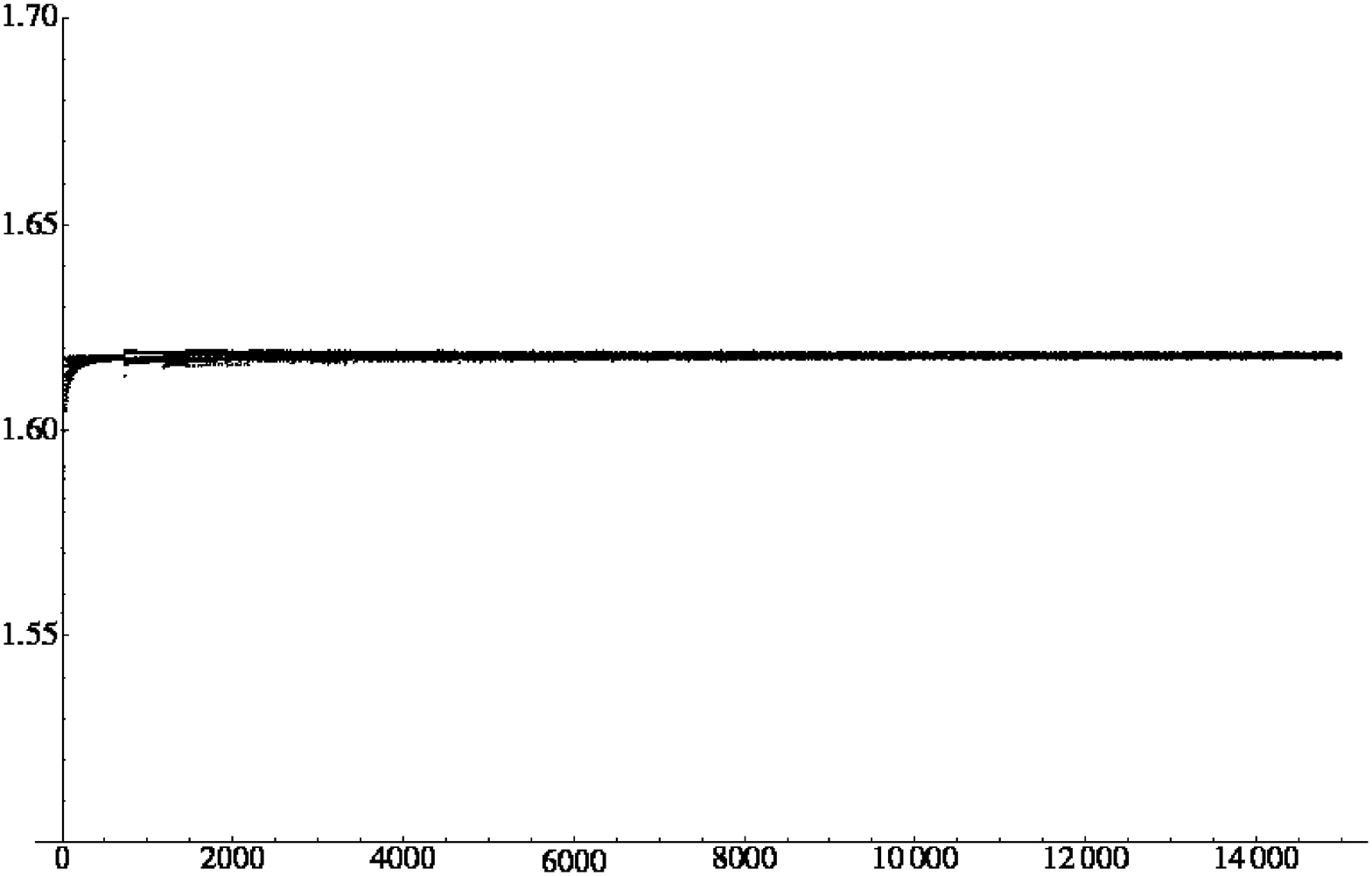}}  
\end{center}
\begin{center}
	  {\includegraphics[width=0.8\textwidth]{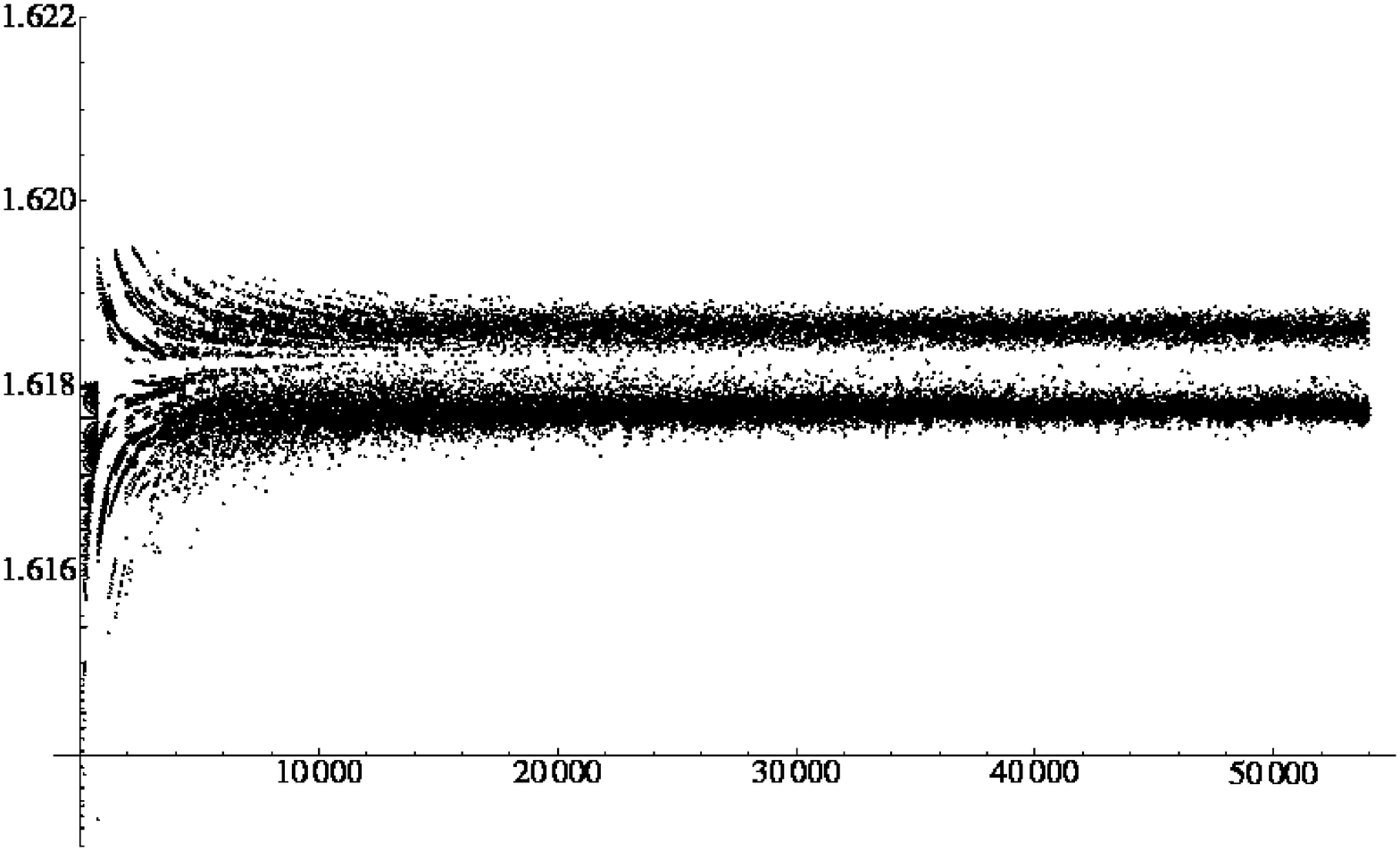}}  
\end{center}\caption{The ratio $b_n/a_n$ for $(731,1183)\G$. We 'expect' to 
see a split since $(731,1183)$ is a splitting pair, indeed, a Wythoff pair. 
Unfortunately, we note that Mathematica has had some problems 
of showing the correct output 
for small $n$ in the lower picture (the upper picture is correct), 
but the splitting tendency for $n$ about 50000 is correctly visualized.} 
\end{figure}

\end{document}